\documentclass[10pt,reqno]{article}
\usepackage{graphicx, latexsym}
\usepackage{amsmath}
\usepackage{amsthm}
\usepackage{amssymb, amsfonts}
\usepackage[numbers]{natbib}
\usepackage{amscd}
\usepackage{empheq}
\usepackage{enumerate}
\usepackage{hyperref}
\usepackage{calc}
\usepackage{accents}
\usepackage{ulem}
\usepackage{tikz}
\usepackage{epstopdf}
\usepackage{caption}
\usepackage{subcaption}
\usetikzlibrary{shapes.geometric, arrows}
\usepackage[toc,page]{appendix}

\graphicspath{ {Figs/} }

\makeatother
\graphicspath{ {Figs/} }

\newtheorem{theorem}{Theorem}[section]

\newtheorem{Remark}{Remark}[section]
\newtheorem{lemma}{Lemma}[section]

\numberwithin{equation}{section}

\newcommand{\R}{\mathbb{R}}
\newcommand{\ms}{\mathrm{s}}
\newcommand{\mY}{\mathrm{Y}}
\newcommand{\mv}{\mathrm{v}}

\newcommand{\bM}{{\bf M}}

\newcommand{\bX}{{\bf X}}
\newcommand{\bw}{{\bf w}}
\newcommand{\bD}{{\bf D}}
\newcommand{\bJ}{{\bf J}}
\newcommand{\bP}{{\bf P}}

\newcommand{\bV}{{\bf V}}
\newcommand{\bH}{{\bf H}}
\newcommand{\bL}{{\bf L}}

\newcommand{\f}{{\bf f}}
\newcommand{\bu}{{\bf u}}

\newcommand{\bn}{{\bf n}}
\newcommand{\bU}{{\bf U}}

\newcommand{\bK}{{\bf K}}
\newcommand{\bz}{{\bf z}}

\newcommand{\bO}{{\bf O}}

\newcommand{\btau}{\mbox{\boldmath $\tau$}}

\newcommand{\brho}{\mbox{\boldmath $\rho$}}

\newcommand{\bzeta}{\mbox{\boldmath $\zeta$}}
\newcommand{\bchi}{\mbox{\boldmath $\chi$}}
\newcommand{\bLambda}{\mbox{\boldmath $\Lambda$}}

\newcommand{\U}{\mbox{\boldmath ${\mathcal{V}}$}}                                                                
\newcommand{\bPi}{\mbox{\boldmath $\Pi$}}
\newcommand{\bnu}{\mbox{\boldmath $\nu$}}

\newcommand{\bv}{\mbox{\boldmath ${\bw}$}}                                                                
\renewcommand{\psi}{\mbox{${\theta}$}}

\newcommand{\e}{{\bf e}}

\DeclareMathOperator*{\esssup}{ess\,sup}
\DeclareMathOperator*{\maximum}{max}

\newcommand{\vertiii}[1]{{\left\vert\kern-0.25ex\left\vert\kern-0.25ex\left\vert #1 
    \right\vert\kern-0.25ex\right\vert\kern-0.25ex\right\vert}}

\begin{document}
\thispagestyle{plain}

\title
{ \large\bf A discontinuous Galerkin method for the evolutionary Navier-Stokes-Darcy model}

\author{Sweta Chakraborty\thanks{Department of Mathematical Sciences, Tezpur University, Tezpur, Sonitpur, Assam-784028, India. Email: msp21103@tezu.ac.in}, ~~
Deepjyoti Goswami\thanks{Department of Mathematical Sciences, Tezpur University, Tezpur, Sonitpur, Assam-784028, India. Email: deepjyoti@tezu.ernet.in}, ~~and~~
Saumya Bajpai\thanks{School of Mathematics and Computer Science, Indian Institute of Technology Goa, Ponda, Goa-403401, India. Email: saumya@iitgoa.ac.in}}

\date{}
\maketitle

\textbf{Abstract:} In this paper, we have analysed a discontinuous Galerkin (DG) method for the coupled surface-subsurface flow, which is modelled by the evolutionary Navier-Stokes-Darcy (NSD) equations. While the DG method has been applied for the spatial discretization of both Navier-Stokes and Darcy equations, for time, we have employed the first-order backward Euler scheme. The fully discrete error analysis has been carried out for both velocity and pressure of the Navier-Stokes equations, and for the piezometric head of the Darcy equation, and optimal estimates have been derived. Our analysis relies on appropriate treatments of the nonlinear convective term from the Navier-Stokes equations, along with the interface conditions that couple the system. 
We have employed a modified Oseen operator on an appropriate DG setup along with the standard duality argument to derive optimal error estimates in the $L^2$-norm for the unknowns, which are novel for the NSD equations within the DG framework and are supported by numerical simulations demonstrating the scheme's efficiency.
\vspace{0.25cm}


\section{Introduction}\label{s1}

\noindent
The combination of surface and subsurface flows, or incompressible flows over porous media, has garnered a lot of interest among the research community due to its numerous applications in environmental sciences, health sciences, bioengineering, hydrology, and industrial applications (see \cite{D04, LSY02} and references therein). In order to simulate such physical phenomena, the Navier-Stokes equations (NSEs) (or Stokes equations) in the free flow 
region are coupled with the Darcy equations in the porous media. This coupling is accomplished through the interface conditions, such as continuity of the normal velocity, balance of normal forces, and the Beavers-Joseph-Saffman (BJS) law (\cite{S71,LS75}) for the tangential velocity.

\noindent
Both the coupled models, namely, Stokes-Darcy (SD) and Navier-Stokes-Darcy (NSD), have been studied on several occasions, involving the well-posedness and continuous Galerkin (CG) finite element approximations. We refer here to non-exhaustive lists \cite{BC09,CGHHWZ10,D05,DMQ02,LSY02} and \cite{BDQ10,CGR13,CDG24,D04,DQ09} for the SD and NSD models, respectively, and references therein. But when it comes to DG approximations for these models, one finds the literature to be limited.

\noindent
DG methods offer several advantages over CG methods, including flexibility in the design of complex geometries, local conservation of mass, higher-order local approximations, $hp$-adaptivity, and so on, and it becomes quite relevant to adapt this numerical method to analyse our coupled model problem. We find that most of the studies related to the coupled models revolve around the CG method, and very few can be attributed to the study in the DG setup. 
In \cite{ChR09}, the DG methods have only been adopted for the porous media, while CG has been applied in the incompressible flow region. The rigorous mathematical analysis of DG methods for both domains in the case of steady NSD has been done in \cite{GR09}. There, under appropriate smallness restrictions on the given data, the authors have proved the existence of a weak solution and established {\it a priori} estimates. Also, local uniqueness has been established when the data fulfil additional smallness conditions. Then, using the DG approximation for the steady coupled NSD model, error bounds in the energy norm have been derived. For the time-dependent NSEs and the steady Darcy equations, primal DG methods are considered in \cite{CR09, CGPR17} along with Crank-Nicolson and backward Euler schemes for time discretization, respectively.
The existence and uniqueness of discrete solutions have been proved under a suitable smallness condition on the given data. And {\it a priori} error estimates have been derived, which are optimal only in the energy norm. Recently, a strongly conservative and pressure-robust hybridizable DG method for transient and steady NSD problems has been considered in \cite{CLR24} and  \cite{CR23}, respectively.

\noindent
We observe here that, in the DG literature for both steady and unsteady NSD problems, the error analysis has been carried out only in energy norms, whereas numerical experiments show optimal convergence rates in both energy and $L^2$-norms. 
Indeed, as far as we know, the optimal error bounds of the discrete velocity and pressure of the NSEs and the discrete piezometric head of the Darcy in $L^2$-norms are absent from the DG literature for both steady and unsteady NSD problems. Moreover, in the unsteady case, optimal error estimates in the $L^\infty(X)$-norms are missing as well, where the $X$-norm represents the appropriate broken energy norm; for example, see (\ref{Norm1}).

\noindent
In this work, we have considered the fully evolutionary NSD problem, where both the NSEs and the Darcy equations are time-dependent, and to the best of our knowledge, it is taken up for the first time to study in the DG setup. We have established here new and optimal error bounds in both energy and $L^2$-norms, and as such, our investigation of the DG method for the fully evolutionary model can be considered as a first attempt in that direction. 

\noindent Our DG scheme for the evolutionary NSD model follows the one presented in \cite{CGPR17}. The consistency and stability of the scheme follow from there as well. We note that due to the interface conditions of the NSD model, the nonlinear convection of the NSEs in the energy balance is not completely eliminated (see (\ref{nonliearity1})). This necessitates small data conditions; for establishing both {\it a priori} and error estimates, see \cite{CGPR17}. 

\noindent 
Our error analysis focuses on careful handling of the nonlinear terms and the interface terms as they play crucial roles in obtaining optimal estimates. To illustrate our point, let us consider a boundary integral $I=\int_e (u-\Pi_h u)v\,ds$, for an edge $e\in\Gamma$; $\Gamma$ is the interface, and for an $L^2$ projection $\Pi_h$. It can be estimated in the following suboptimal ways: For $u\in H^{k+1},\, v\in H^1$,
\begin{align*}
|I|&\leq \sum_{e\in \Gamma}\|u-\Pi_h u\|_{L^2(e)}\|v\|_{L^2(e)}\leq Ch^{k+1/2}\|u\|_{k+1}\|v \|_1,\\
|I|&\leq \sum_{e\in \Gamma}|e|^{1/4}\|u-\Pi_h u\|_{L^2(e)}\|\|v\|_{L^4(e)}\leq Ch^{k+3/4}\|u\|_{k+1}\|v\|_1,\\
|I|&\leq \sum_{e\in \Gamma}|e|^{1/2}\|u-\Pi_h u\|_{L^4(e)}\|v\|_{L^4(e)}\leq Ch^{k+3/4}\|u\|_{k+1}\|v\|_1.
\end{align*}
In standard DG analysis, boundary integrals are handled by absorbing one of the terms in the energy norm, giving rise to an extra $|e|^{1/2}$ that compensates for the suboptimal results. But in our case, the absence of an interface term in the definition of the energy norm (see (\ref{Norm1})), robs us of this idea, and we are stuck with at most $O(h^{3/4})$.

\noindent
As for the nonlinear convective term of the NSEs, Green's theorem is applied in order to derive optimal $\bL^2$-error bounds; see the DG and DG-two grid literature for NSE, \cite{BGR23,RGB24}. But regarding the NSD, such a mechanism is absent in the literature. And an application of Green's theorem to the nonlinear term will simply introduce an interface term, which would produce suboptimal error bounds (see Remark \ref{bypartsresult}). 

\noindent
Here, our primary focus lies on overcoming these difficulties and recovering optimal results in the fully discrete case. We have begun by introducing an $\bL^2$-projection $\bP_h$ onto a suitable DG setup, utilizing the approximation operator $\bPi_h$, and have proved appropriate approximation properties. Next, following the idea of the linearized dual problem in \cite{CMX09}, 
we have constructed modified Oseen operators $\bO_{fh}$ for velocity and $O_{ph}$ for piezometric head, respectively, on appropriate DG spaces. These operators, $\bO_{fh}$ and $O_{ph}$, play a crucial role in handling the above-mentioned difficulties, thereby recovering optimal convergence rates (see Remark \ref{Oseenhelp}). We have established optimal approximation properties for both $\bO_{fh}$ and $O_{ph}$ by means of standard duality arguments. We would like to point out here that although a linearized dual problem is standard in the DG analysis for nonlinear problems like NSEs, what separates our case is the appropriate choice of the interface terms, leading to optimal error estimates.

\noindent
Then, temporal discretization has been performed using the backward Euler scheme, and fully discrete optimal error $L^\infty(\textbf{L}^2)$ estimates for velocity and $L^\infty({L}^2)$ estimates for piezometric head have been shown. Additionally, optimal error estimates in the energy norm have been derived as well.
Finally, we have established an optimal $L^{\infty}(L^2)$-pressure estimate for the NSEs. These estimates are technically more involved (see Remark \ref{pressure}), and we have provided comments on these effects in Section $5.2$, see Lemmas \ref{Intermidiate1} and \ref{Intermidiate2}. We would like to point out here that our analysis has been carried out for the two-dimensional case, but a similar analysis will go through for the three-dimensional case as well; see Remark \ref{3D}. We end our discussion by simulating two problems, 
verifying the optimal convergence rates, and a benchmark problem,  ensuring the efficiency of our numerical scheme.

\noindent Below, we present our primary contributions:
\begin{itemize}
\item Modified Oseen operators $\bO_{fh}$ for velocity and $O_{ph}$ for piezometric head, and their optimal approximation properties.
\item Optimal $L^\infty(\textbf{L}^2)$ and $L^\infty({L}^2)$-norms, and energy norm error bounds for the fully discrete DG approximations to the velocity and piezometric head, respectively.
\item $L^\infty({L}^2)$-norm error bound for fully discrete DG approximations to the pressure.
\item Numerical experiments validating the theoretical results and efficiency of the scheme.
\end{itemize}
\noindent The rest of the paper is structured as follows: In Section \ref{s2.}, the model problem has been described along with notations, variational formulation, basic assumptions, and standard Sobolev inequalities. Numerical scheme has been provided in Section \ref{s3.}. Section \ref{s4.} contains the properties of approximation operators, while the fully discrete optimal error analysis for velocity, piezometric head, and pressure has been dealt with in Section \ref{s6.}. The numerical results have been discussed in Section \ref{s7.}. Finally, conclusions have been drawn in Section \ref{s8.}.

\noindent
\section{Preliminaries}\label{s2.}
\subsection{Model Problem}
Consider a coupled fully evolutionary NSD model in a bounded connected domain $\Omega=\Omega_f\cup \Omega_p$ in $ \mathbb{R}^d (d=2,3)$, where $\Omega_f$  and $\Omega_p$ represent the fluid flow region and the porous media region, respectively, satisfying $\Omega_f\cap\Omega_p=\emptyset$. We further assume that each subdomain is convex.
 Let $\partial\Omega$, $\partial\Omega_f$, and $\partial\Omega_p$ denote the boundaries of the domain $\Omega$, and subdomains $\Omega_f$ and $\Omega_p$, respectively. The interface between the two regions is given by ${\Omega}_f\cap{\Omega}_p=\Gamma$. Finally, the exterior boundaries of the free flow region and porous medium region are given by $\Gamma_f=\partial\Omega_f\cap\partial\Omega$ and $\Gamma_p=\partial\Omega_p\cap \partial\Omega$, respectively. 
\begin{figure}[h!]
    \centering
    \scalebox{0.60}{
    \begin{tikzpicture}
    \draw (0,0) rectangle (10,5);
 \draw [color={rgb,255:red,38; green,7; blue,211}, draw opacity=1, thick]
              (0,2.5) .. controls (4,4) and (6,1) .. (10,2.5)
              node[sloped,inner sep=0cm,above,pos=.235,
              anchor=south west,
              minimum height=(2.0)*0.3cm,minimum width=(3.0)*.3cm](N){};
 \coordinate (P1) at (4.7,2.5);
        \coordinate (P2) at (5.3,2.5);
\draw[->, thick, black] (P2) -- ++(0,-0.6) node[left] {$\bn_{\Gamma}$};
 \draw[->, thick, black] (P2) -- ++(0.7,0) node[above] {$\btau_{\Gamma}$};
 \node at (5, -0.5) {$\Gamma_{p}$};
        \node at (-0.5, 1.25) {$\Gamma_{p}$};
        \node at (10.5, 1.25) {$\Gamma_{p}$};
        \node at (5, 5.5) {$\Gamma_{f}$};
        \node at (-0.5, 4.25) {$\Gamma_{f}$};
        \node at (10.5, 4.25) {$\Gamma_{f}$};
        \node at (6.8, 2.5) {$\Gamma$};
        \node at (4, 1.5) {$\Omega_p$};
        \node at (6.5, 4) {$\Omega_f$};    
\end{tikzpicture}}
    \caption{Coupled domain of fluid flow and porous media flow.} 
    \label{domain}
\end{figure}
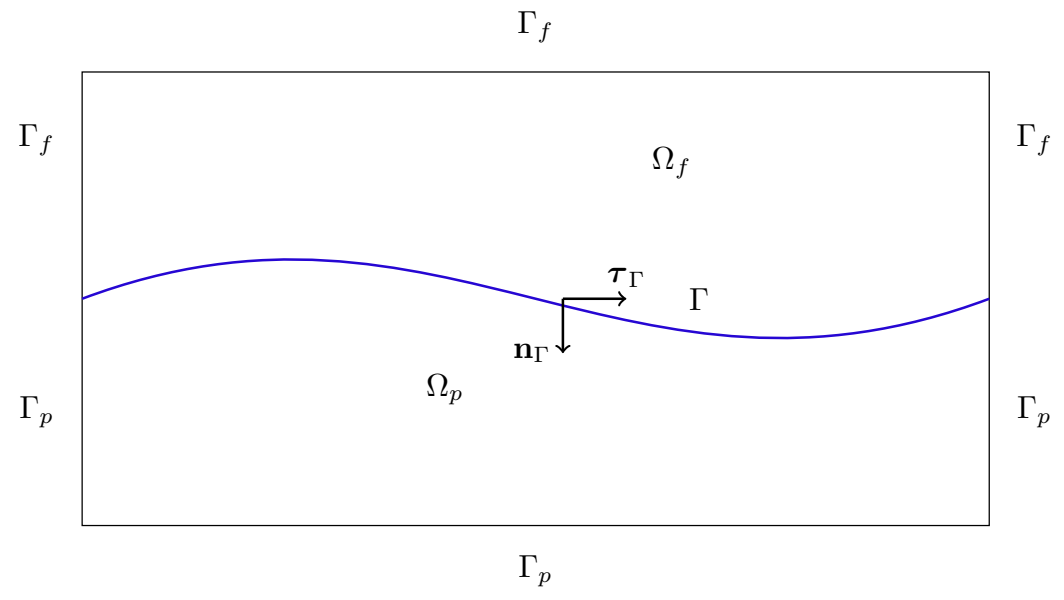

\noindent
The Navier-Stokes equations, which govern the fluid flow in $\Omega_f$, are expressed as \vspace*{-2mm}
\begin{equation}
    \left. %
\begin{aligned}
     \bu_t- \nabla\cdot(2\mu \textbf{D}(\bu)-p{\bf I})+ \bu \cdot \nabla \bu &= \f_f ~~\mbox{in}~ \Omega_f,\,t>0,\\
\nabla \cdot \bu &= 0~~\mbox{in}~  \Omega_f,\,t>0, \\
\bu(0)& = \bu_0~~\mbox{in}~ \Omega_f, \label{NSEDARCY1}
\end{aligned}\right\}%
\end{equation}
\noindent where $\bu_t$ denotes differentiation of $\bu$ in time. Here, $\bu=(u_i)_{i=1}^d$ indicates the velocity vector, $p$ the pressure of the fluid in $\Omega_f$, coefficient $\mu>0$ the kinetic viscosity, $\textbf{D}(\bu)=(1/2)(\nabla\bu+(\nabla\bu)^T)$ the deformation tensor with ${\bf I}$ being the identity tensor, and $T$ is the final time with $0<T<\infty$. The forcing term $\f_f$ is a given function in the domain $\Omega_f$. Other notations have the usual meaning as expressed below: \vspace*{-1mm}
\[\bu\cdot\nabla\bu=\sum_{i=1}^d u_i\frac{\partial\bu}{\partial x_i}\quad \text{and}\quad \nabla\cdot \bu=\sum_{i=1}^d \frac{\partial u_i}{\partial x_i}.\]
The flow in porous media is characterized by the unsteady Darcy equations, which are given by \vspace*{-2mm}
\begin{equation}
    \left. %
\begin{aligned}
 S_0\phi_t-\nabla \cdot ( \bK\nabla \phi) & = f_p ~~\text{in}~~\Omega_p,\,t>0,\\
	\phi(0) &= \phi_0~~ \mbox{in}~~ \Omega_p,\label{NSEDARCY7}
	\end{aligned}\right\}%
\end{equation}
\noindent
where $\phi$ is the piezometric head of the porous media flow in $\Omega_p,\, S_0>0$ is the specific mass storativity coefficient, and $f_p$ is the forcing term in $\Omega_p$. $\bK\in\R^{d\times d}$ is the permeability tensor, which is presumed to be symmetric positive definite, with eigenvalues uniformly bounded above and away from zero.
Without loss of generality, we impose homogeneous Dirichlet boundary conditions on $\Gamma_f$ and $\Gamma_p$: \vspace*{-1mm}
\begin{align}
	\bu ={\bf 0} \quad \mbox{on }\ \Gamma_f,\,t\geq 0,\quad
	 \phi =0 \quad \mbox{on }\ {\Gamma_{p}},\,t\geq 0. \label{NSEDARCY9}
	 	\end{align}
The two flows in each subdomain are coupled by the following conditions across the interface $\Gamma$: \vspace*{-1mm}
\begin{align}
\bu\cdot\bn_{\gamma}  = -\bK\nabla \phi\cdot\bn_{\gamma} \quad &\mbox{on } \Gamma,\,t>0,\label{INT1} \\
((-2\mu \textbf{D}(\bu)+p{\bf I})\bn_{\gamma})\cdot\bn_{\gamma} =  g\phi\quad&\mbox{on }\Gamma,\,t>0,\label{INT2} \\
 \bu \cdot{\btau}^{j}_{\gamma} = -2\mu G^j(\textbf{D}(\bu)\bn_{\gamma})\cdot\btau^j_{\gamma},\quad 1\leq j\leq d-1\quad &\mbox{on }~\Gamma,\,t>0, \label{INT3}
\end{align}
where $G^j=\frac{\mu\alpha}{(\bK\btau^j_{\gamma},\btau^j_{\gamma})^{1/2}}$, $\alpha>0$ is determined from experimental data, $\bn_{\gamma}$ is the unit outward normal vector to $\Gamma$ pointing from $\Omega_f$ to $\Omega_p$, $\{\btau^j_{\gamma}\}_{j=1}^{d-1}$ are linearly independent unit tangential vectors on $\Gamma$, and $g$ is the gravitational acceleration. The condition (\ref{INT1}) presents the conservation of mass across the interface ${\Gamma}$. The second condition (\ref{INT2}) guarantees the balance of forces across the interface, and  the third one (\ref{INT3}) is known as the BJS condition, which depicts that the tangential component of the velocity is proportional to the shear stress. For more details of the fully evolutionary model, we refer to \cite{LD22,XH20}.
\subsection{Functional Settings} 
In the entire paper, for a domain $\mathfrak{D}\subset\mathbb{R}^d$, boldface characters are used to indicate $\mathbb{R}^d$-vector-valued functional spaces such as $\bH^1(\mathfrak{D})=(H^1(\mathfrak{D}))^d$ and $\bL^2(\mathfrak{D})=(L^2(\mathfrak{D}))^d$. We consider, for $1\le r\le\infty$, the Lebesgue integrable spaces $L^r(\mathfrak{D})$ with norm $\|\cdot\|_{L^r(\mathfrak{D})}$ and the classical Sobolev spaces $W^{m,r}(\mathfrak{D})$, for any integer $m\geq 0$ \cite{A75}:
\[W^{m,r}(\mathfrak{D})=\{\mv\in L^r(\mathfrak{D}):\, \partial^\ms \mv\in L^r(\mathfrak{D}) \,\,\forall\, |\ms|\leq m\}, \] 
\noindent where $\ms$ represents a multi-index $(\ms_i)_{i=1}^d$ with $|\ms|=\sum \ms_i$ and $\partial^\ms \mv$ the partial derivatives of $\mv$ of order $|\ms|$. The associated norm on $W^{m,r}(\mathfrak{D})$ is denoted by $\|\cdot\|_{m,r,\mathfrak{D}}$, and the semi-norm by $|\cdot|_{m,r,\mathfrak{D}}$. For $r=2$, $W^{m,2}(\mathfrak{D})$ is a Hilbert space $H^m(\mathfrak{D})$, and $r$ is dropped from the subscripts of its norm and semi-norm. For time-dependent measurable functions $\mv$ with $\mv(t)$ in a Banach space $\mY$ for $t\in[0,T]$, we consider the space $L^r(0,T;\mY)$; $\mv\in L^r(0,T;\mY)$ satisfying
\[\int_{0}^{T} \|\mv\|_{\mY}^{r} dt< \infty ~\mbox{for}~ 1< r< \infty ~~\mbox{and}~~\esssup_{0\leq t\leq T } \|\mv(t)\|_{\mY}<\infty \,\mbox{ for }\, r=\infty.\] 
The required functional spaces for the weak solution $(\bu,p,\phi)$ of (\ref{NSEDARCY1})-(\ref{INT3}) are given by \vspace*{-2mm}
\begin{align*}
\bX &= \{\bv\in (H^1(\Omega_f))^d: \bv = 0 ~\text{on}~{\Gamma}_f\},\quad
M_f = L_0^2(\Omega_f),\quad
 M_p = \{\psi\in H^1(\Omega_p): \psi = 0 ~\text{on}~\Gamma_p\}.
\end{align*}%
\noindent We end this segment by recollecting standard inequalities such as Poincar{\`e}'s, Korn's, and Sobolev inequalities \cite[pp. 2057:($1.21$)-($1.23$)]{GR09}: there exist constants $d_1,d_2,d_3>0$ that depend only on $\Omega_f$, such that $\forall{\bz}\in \bX$,\vspace*{-2mm}
\[| \bz|_{H^1(\Omega_f)}  \leq d_1\|\textbf{D}(\bz)\|_{L^2(\Omega_f)},\:
\| \bz\|_{L^p(\Omega_f)}  \leq d_2|\bz|_{H^1(\Omega_f)},(p=2,3,4), \: \| \bz\|_{L^r(\Gamma)}\leq d_3 |\bz|_{H^1(\Omega_f)},\,(r=2,4).\]%

\section{Numerical Scheme}\label{s3.}
\noindent In order to discretize the coupled problem, we presume a family of triangulations $\mathcal{E}_{h}$ of $\bar{\Omega}$ which is regular and comprises $E$'s, triangles if $d=2$ or tetrahedra if $d=3$ of diameter $h_E$. Let $h=\maximum_{E\in \mathcal{E}_{h}}h_E$. We consider that the interface $\Gamma$ consists of a finite union of the triangle edges. Consequently, the restriction of $\mathcal{E}_h$ to $\Omega_i$, denoted by $\mathcal{E}_{ih}$, for $i=f,p$, is also considered as a regular family of triangulations of $\bar{\Omega}_i$. 
Regular indicates that there exists a constant $\beta>0$ which does not depend on $h$, such that (see \cite{C02}),
\begin{equation*}
\frac{h_E}{\varrho_E}\leq \beta, \quad \forall E\in \mathcal{E}_h,
\end{equation*}
\noindent where $\varrho_E$ represents the diameter of the circle which is inscribed in $E$.
By $\Gamma_{ih}$, we denote the set of all edges of $\mathcal{E}_{ih}$ that are interior to $\Omega_i$ for $i=f,p$. With each interior edge $e$, we assign a unit normal vector $\bn_e$. For any boundary edge
$e\in \Gamma_f\cup\Gamma_p$, the associated unit vector $\bn_e$ coincides with the outward normal vector. For $e\in\Gamma$, we choose $\bn_e=\bn_{\gamma}$ such that $\bn_e$ is the unit outward normal pointing from $\Omega_f$ to $\Omega_p$. When a unit normal vector $\bn_e$ is directed from an element $E_1$ to an element $E_2$, then we describe the jump $[\omega]$ and average $\{\omega\}$ of a function $\omega$ as
\[ [\omega]_e=(\omega|_{E_1})|_e-(\omega|_{E_2})|_e, \quad \{\omega\}_e=\frac{1}{2}(\omega|_{E_1})|_e+\frac{1}{2}(\omega|_{E_2})|_e.\] 
For an edge $e\in\Gamma_f\cup\Gamma_p$, we have $[\omega]_e=\{\omega\}_e=\omega|_e$.  

\noindent We next introduce the following finite dimensional discrete subspaces $\bX_h$, $M_{fh}$, and $M_{ph}$:
\begin{eqnarray*}
 \bX_h &=& \{\bv_h\in \bL^2(\Omega_f):\, \bv_h|_{E}\in \mathbb{P}_{k_1}(E)^2,\, \forall E\in \mathcal{E}_{fh}\},\\
 M_{fh} &=& \{q_h\in L_0^2(\Omega_f):\, q_h|_{E}{\in \mathbb{P}_{k_1-1}}(E),\, \forall E\in \mathcal{E}_{fh}\},\\
 M_{ph} &=& \{\psi_h\in L^2(\Omega_p):\, \psi_h|_{E}{\in \mathbb{P}_{k_2}}(E),\, \forall E\in \mathcal{E}_{ph}\},
\end{eqnarray*}
\noindent where $\mathbb{P}_k(E)$ presents the polynomial spaces of degree $k$ on $E$ with two positive integers $k_1\geq 1,k_2\geq 1$.

\noindent For ${M_{fh}}$, we take the $L^2$ induced norm, and for the other spaces, we associate the following norms:
\begin{align}
\|\bv_h\|^2_{\bX_{1}}&=\bigg(2\sum_{E\in\mathcal{E}_{fh}}\|\textbf{D} (\bv_h)\|^2_{L^2(E)}+\sum_{e\in \Gamma_{fh}\cup\Gamma_f} \frac{\sigma_e}{|e|}\|[\bv_h]\|^2_{L^2(e)}\bigg),\: \forall \,\bv_h\in \bX_{h},\label{Norm1} \\
\|\psi\|^2_{M_{1p}}&=\bigg(\sum_{E\in\mathcal{E}_{p}}\|\bK^{1/2} \nabla \psi_h\|^2_{L^2(E)}+\sum_{e\in \Gamma_{ph}\cup\Gamma_p} \frac{\sigma_e}{|e|}\|[\psi_h]\|^2_{L^2(e)}\bigg),\: \forall \psi\in {M_{ph}}. \label{Norm2}
\end{align}
\noindent Here, $|e|$ indicates the length of edge $e$. The penalty parameter $\sigma_e>0$ depends on $e$ but is independent of $h$.
We further present the discrete bilinear forms corresponding to the DG formulation as follows: $a_{NS}:\bX_h\times\bX_h\rightarrow \mathbb{R}$,  $b_{NS}:\bX_h\times M_h\rightarrow \mathbb{R}$ and $a_{D}:{M_h}\times {M_h}\rightarrow \mathbb{R}$: 
\begin{align}
\quad a_{NS}(\bu_h,\bv_h)& = 2\mu\sum_{E\in \mathcal{E}_{fh}}( \textbf{D}(\bu_h), \textbf{D}(\bv_h))_{E} -2\mu\sum_{e\in \Gamma_{fh}\cup\Gamma_f} (\{ \textbf{D}(\bu_h)\bn_e\}, [\bv_h])_e\nonumber \\
& +2\epsilon_1\mu\sum_{e\in \Gamma_{fh}\cup\Gamma_f}( \{\textbf{D}(\bv_h)\bn_e\}, [\bu_h] )_e
 + \mu \sum_{e\in \Gamma_{fh}\cup\Gamma_f} \frac{\sigma_e}{|e|}([\bu_h],[\bv_h])_e\nonumber,\\
 b_{NS}(\bv_h,q_h) &= -\sum_{E\in \mathcal{E}_{fh}}( q_h,\nabla \cdot \bv_h)_e +\sum_{e\in \Gamma_{fh} \cup \Gamma_f}( \{q_h\},[\bv_h]\cdot \bn_e)_e,\nonumber \\
 a_{D}(\phi_h,\psi_h) &= \sum_{E\in \mathcal{E}_{ph}}( \bK \nabla \phi_h, \nabla\psi_h)_E -\sum_{e\in \Gamma_{ph}\cup\Gamma_p} (\{\bK \nabla\phi_h\cdot\bn_e\}, [\psi_h] )_e\nonumber \\
& +\epsilon_2\sum_{e\in \Gamma_{ph}\cup\Gamma_p}( \{ \bK\nabla \psi_h\cdot\bn_e\}, [\phi_h])_e+  \sum_{e\in \Gamma_{ph}\cup\Gamma_p} \frac{\sigma_e}{|e|}([\phi_h],[\psi_h])_e,\nonumber 
\end{align}
\noindent where for $i=1,2$, $\epsilon_i\in\{-1,0,1\}$. For $\epsilon_1=\epsilon_2=1$, it becomes the non-symmetric interior penalty Galerkin (NIPG) method. In case the parameters $\epsilon_i,~ i=1,2$ take the value $-1$ (or $0$), the method is called the symmetric (or incomplete) interior penalty Galerkin method, and in short, SIPG (or IIPG) method.
 As for the trilinear form $c_{NS}(\cdot,\cdot;\cdot,\cdot)$ for the convection term, it has been discretized by applying upwinding of Lesaint-Raviart (see \cite{LR74}), introduced in \cite{GRW005}, as
\begin{align}
c_{NS}(\bnu_h,\bu_h;\bz_h,\bw_h) = n(\bu_h;\bz_h,\bw_h)+l(\bnu_h,\bu_h;\bz_h,\bw_h),\label{trilinear}
\end{align}\vspace*{-2mm}
\text{where}
\small{\begin{align}
&n(\bu_h;\bz_h,\bw_h) =  \sum_{E\in\mathcal{E}_{fh}}  (\bu_h\cdot \nabla \bz_h, \bw_h)_E  +\frac{1}{2}\sum_{E\in \mathcal{E}_{fh}}(\nabla\cdot \bu_h,\bz_h\cdot \bw_h)_E -\frac{1}{2}\sum_{e\in \Gamma_{fh}\cup\Gamma_f}([\bu_h]\cdot \bn_e,\{\bz_h\cdot \bw_h\})_e ,
\nonumber \\
&l(\bnu_h,\bu_h;\bz_h,\bw_h)  =  (|\{\bu_h\}\cdot \bn_E|(\bz_h^{int}-\bz_h^{ext}), \bw_h^{int})_{\partial E_-(\bnu_h) \setminus \Gamma}, \,\,\forall \bu_h,\bz_h,\bw_h \in \textbf{X}_{h}.\nonumber
\end{align} }
\noindent Here, $\partial E_-(\bnu_h)$ denotes the inflow boundary of the element $E$ with respect to a function $\bnu_h\in\bX_{h}$ expressed as
\[\partial E_-(\bnu_h)=\{x\in \partial E:\{\bnu_h\}\cdot \bn_E<0\},\]
and $\bz^{int}$ (respectively $\bz^{ext}$) indicates the trace of the function $\bz$ on a side of $E$ coming from the interior of $E$ (respectively, exterior of $E$) at any point $x$ of $\partial E$.

\noindent Before proceeding further, regarding the energy norm, we recollect the following bound for the elements in the broken-$H^1$ space $\bX_1$ $= \{\bv\in \bL^2(\Omega_f):\, \bv|_{E}\in \bH^1(E),\forall E\in\mathcal{E}_{fh}\}$ \cite[Lemma 6.2]{GRW005}:\vspace*{-2mm} 
\begin{align}
\|\bv\|_{L^{\mathfrak{p}}(\Omega_f)}\leq C_{\mathfrak{p}}\|\bv\|_{\bX_1}, \: \forall\bv\in\bX_1, \: \mbox{ for }\, {\mathfrak{p}}\in[2,\infty),\label{Poincare}
\end{align}
\noindent presuming that $C_{\mathfrak{p}}$ is a positive constant.
In the three-dimensional case, the result is valid for $2\leq\mathfrak{p}\leq 6$.

\noindent We now recollect a few approximation properties of the discrete spaces. Also note that henceforth the constant $C>0$ that appears in the subsequent analysis is independent of $h$ and $E$. \begin{lemma}\cite[Lemma 6.1]{CGPR17}\label{Piprojection}
There exists an operator $\bPi_h:\bX\rightarrow\bX_h$ which satisfies \begin{align*}
b_{NS}(\bw-\bPi_h(\bw),q_h)=0 \quad \forall \bw\in \bX, \,\,\forall q_h\in M_{fh},
\end{align*}
and $\forall E  \in {\mathcal{E}}_{fh}, \bw\in \bX\cap W^{s,r}(E)^d,\,\, 1\leq r \leq \infty, \,\, 1\leq s \leq k_1+1$,
\begin{align*}
\|\bw\,-\,\bPi_h(\bw)\|_{L^r(E)} + h_E\|\nabla(\bw\,-\,\bPi_h(\bw))\|_{L^r(E)}&\leq Ch_E^s |\bw|_{W^{s,r}(\Delta_ E)},
\end{align*}
\noindent where a macro-element $\Delta_ E\subset \Omega_f$ is used to construct $\bPi_h(\bw)$ in $E$. Additionally, $\forall s, 1 \leq s \leq k_1 + 1$, $\bPi_h$ satisfies the following bound:
\begin{align}
\|\bw-\bPi_h(\bw)\|_{L^2(\Omega_f)} + h\,\|\bw-\bPi_h(\bw)\|_{\bX_1}&\leq Ch^s |\bw|_{{s},\Omega_f}, \,\,\forall \bw\in\,\bH^s(\Omega_f)\cap\bX.\label{Piproperty}
\end{align}
\end{lemma}%
\noindent Moreover, there exists a local projection operator $r_h:L_0^2(\Omega_f)\rightarrow M_{fh}$ (see \cite[(3.10)-(3.11)]{GRW005}), defined for any $E\in \mathcal{E}_{fh}$ and $z_h\in \mathbb{P}_{r-1}(E)$, \vspace*{-2mm}
\begin{align*}
\int_E z_h{(y-r_h y)}dx=0, \quad\forall y\in L_0^2{(\Omega_f)}.
\end{align*}
\noindent Additionally, $r_h$ satisfies 
\begin{align}
\| y-r_h( y)\|_{L^2(E)}\leq Ch^s_E | y|_{s,E},\quad \forall s\in [0,r],\,\, \forall  y\in H^s(\Omega_f)\cap L_0^2(\Omega_f).\label{projection1}
\end{align}
\noindent To interpolate the piezometric head in Darcy region, we use the Lagrange interpolant $L_h$ in $\mathbb{P}_{k_2}$(see \cite[$(38)$]{CGPR17}) satisfying the following error bounds: for any $E\in \mathcal{E}_{ph}$ and for $r\in [1, \infty),\, s \in  [2,k_2+1]$,
\begin{align}
   \,\,|\psi-L_h(\psi)|_{W^{l,r}( E)}\leq Ch^{s-l} |\psi|_{W^{s,r}( E)},\quad  \forall \,\psi \in  W^{s,r}(E), \; (l=0,1).\label{projection3}
\end{align}
\noindent We apply the backward Euler scheme for temporal discretization. Let us consider a uniform partition $\{t_i\}_{i=0}^N$ in the time interval $[0,T]$. We denote the time step by $\Delta t>0$, and $t_N=T$ is the final time. Therefore, $t_i=i\Delta t$, $0\leq i\leq N$. Further, for a given sequence $\{\nu^i\}_{i\geq 0}$, we define $\partial_t\nu^i$ as follows:
\begin{align*}
\partial_t\nu^i=\frac{1}{\Delta t}(\nu^i-\nu^{i-1}).
\end{align*}
The fully discrete numerical scheme for the NSD equations is presented below: \\
Find $(\bU_h^i,P_h^i,\Phi_h^i)_{i\geq 1}\in(\bX_h,M_{fh},M_{ph})$, such that $\forall(\bv_h,q_h,\psi_h)\in(\bX_h,M_{fh},M_{ph})$,
\small{%
\begin{align}
(\partial_t &\bU_h^i,\bv_h)+ gS_0(\partial_t \Phi_h^i,\psi_h)+ a_{NS}(\bU_h^i,\bv_h)+ ga_D(\Phi_h^i,\psi_h) 
+ c_{NS}(\bU_h^{i-1},  \bU_h^{i-1};\bU_h^{i},\bv_h)
+ b_{NS}(\bv_h,P_h^i)\nonumber\\
\hspace*{1mm}&+ \sum_{j=1}^{d-1}\frac{1}{G^j}(\bU_h^i\cdot\btau^j_{\gamma},\bv_h\cdot\btau^j_{\gamma})_{{\Gamma}}
+   g(\Phi_h^i,\bv_h\cdot{\bn_{\gamma}})_{{\Gamma}} -  g(\bU_h^i\cdot\bn_{\gamma},\psi_h)_{{\Gamma}}=(\textbf{f}_f^i,\bv_h)_{{\Omega}_f}+ g(f_p^i,\psi_h)_{{\Omega}_p},\label{FNSED1}\\
&\quad\quad\quad\quad\quad\quad\quad\quad\quad\quad\quad\quad\quad\quad\quad\quad\quad\quad\quad\quad\quad b_{NS}(\bU_h^i,q_{h})= \,0, \hspace*{4cm}\label{FNSED2}
\end{align}}%
\noindent given $\bU_h^0$ and $\Phi_h^0$, which are suitable approximations of $\bu_0$ and $\phi_0$, respectively.
Here, $\f_f^i$ and $f_p^i$ are the forcing terms at the $i$-th time level in $\Omega_f$ and $\Omega_p$, respectively. To derive the consistency of the scheme, one can refer to \cite[ Proposition $3.1$]{GR09}. For subsequent analysis, we introduce the discrete divergence-free space, \vspace*{-1mm}
\[\bV_h=\{\bv_h\in\bX_h:\forall  q_h \in M_{fh},\,b_{NS}(\bv_h,q_h)=0\}.\]
\noindent Now restricting the problem (\ref{FNSED1})-(\ref{FNSED2}) in $\textbf{V}_h$, we get the equivalent DG formulation for the coupled system as:
For $\bv_h\in \bV_h$, $\psi_h\in M_{ph}$, we find $\{\bU_h^i\}_{i\geq 1}\in \bV_h$, $\{\Phi_h^i\}_{i\geq 1}\in M_{ph}$ such that 
\small{
\begin{align}
(\partial_t \bU_h^i,\bv_h)+ gS_0(\partial_t \Phi_h^i,\psi_h)+ a_{NS}(\bU_h^i,\bv_h)+ g&a_D(\Phi_h^i,\psi_h) 
+ \sum_{j=1}^{d-1}\frac{1}{G^j}(\bU^i_h\cdot\btau^j_{\gamma},\bv_h\cdot\btau^j_{\gamma})_{{\Gamma}}\nonumber\\
+ c_{NS}(\bU_h^{i-1},  \bU_h^{i-1};\bU_h^{i},\bv_h)
+   g(\Phi_h^i,\bv_h\cdot{\bn_{\gamma}})_{{\Gamma}} & -  g(\bU_h^i\cdot\bn_{\gamma},\psi_h)_{{\Gamma}}=\,(\textbf{f}_f^i,\bv_h)_{{\Omega}_f}+ g(f_p^i,\psi_h)_{{\Omega}_p} .\label{FNSED3}
\end{align}}%
\noindent We next recall the following standard trace inequalities 
to carry out our analysis.
\begin{lemma}\cite[pp. $23$, $(2.1)$-$(2.2)$]{R08},\cite[pp. $118$, $(7.7)$]{KJ98}\label{TL1}
For each $E\in \mathcal{E}_{fh}$, the following inequalities hold: \begin{align}
&\|\nabla^i\bv\|_{0,e}  \leq C(h_E^{-1/2}\|\nabla^i\bv\|_{0,E}+h_E^{1/2}\|\nabla^{i+1} \bv\|_{0,E})\quad \forall e\in \partial E,\,\,\forall \bv\in \bH^{i+1}(E),\,i=0,1,\label{T1}\\
&\|\bv\|_{L^4(e)}  \leq C(h_E^{-1/4}\|\bv\|_{0,E}+h_E^{ (3-d)/4}\|\nabla \bv\|_{0,E})\quad \forall e\in \partial E,\,\, \forall\bv\in \bH^1(E).\label{T3}
\end{align}%
\end{lemma}
\noindent Similar to the trace inequalities defined in (\ref{T1})-(\ref{T3}), we have the following trace and inverse inequalities in the discrete setup as well. \vspace*{-2mm}
\begin{lemma} \label{TL2} \cite{R08}[pp. 23, $(2.3)$-$(2.5)$],\cite[section 1.4.3]{DDE11} For each $E\in \mathcal{E}_{fh}$ and $1\leq p \leq\infty $, the following inequalities hold:\vspace*{-5mm}
\begin{align}
&\|\nabla^i\bv_h\|_{L^2(e)}\leq Ch_E^{-1/2}\|\nabla^i\bv_h\|_{L^2(E)}\quad \forall e\in \partial E, \,\,\forall \bv_h\in\bX_h,\,i=0,1,\label{T4}\\
&\|\nabla\bv_h\|_{L^2(E)}\leq Ch_E^{-1}\|\bv_h\|_{L^2(E)}\quad \forall \bv_h\in\bX_h,\label{T6}\\
&\|\bv_h\|_{L^p(E)}\leq Ch_E^{ d(1/p-1/2)}\|\bv_h\|_{L^2(E)}\quad \forall \bv_h\in\bX_h.\label{T7}
\end{align}
\end{lemma}
\noindent The above-mentioned inequalities hold similarly for the discrete space $\bM_{ph}$. The following estimates yield the 
usual Sobolev embeddings and trace inequalities to 
the discrete spaces (see \cite[(34)-(38)]{CGPR17}): 
\begin{align}
&\text{For}\,\, 2\leq r\leq 6,\,\, \|\bv_h\|_{L^r(\Omega_f)} \leq C  \|\bv_h\|_{\bX_{1}},\quad \|\psi_h\|_{L^r(\Omega_p)} \leq C  \|\psi_h\|_{M_{1p}}, ~\forall \bv_h\in \bX_{h},\,\forall\psi_h\in M_{ph}. \label{T8}\\
&\text{For}\,\, 2\leq r\leq 4,\,\, \|\bv_h\|_{L^r(\Gamma)} \leq C  \|\bv_h\|_{\bX_{1}},\,\,\,\quad\|\psi_h\|_{L^r({\Gamma})} \leq C  \|\psi_h\|_{M_{1p}}, ~\forall \bv_h\in \bX_{h},\, \forall\psi_h\in M_{ph}.\label{T9}
\end{align}
\noindent For $r=2$, the constants of (\ref{T8}) will be represented by $C_2$ and $C_1$, respectively.
Also, with the help of (\ref{T9}), the following term can be bounded:
\begin{align*}
 \frac{1}{2}(\bv_h\cdot\bn_{\gamma},\bz_h\cdot\bz_h)_{\Gamma}\leq C_0 \|\bv_h\|_{\bX_1} \|\bz_h\|^2_{\bX_1}\quad \forall \bz_h, \bv_h\in {\bX}_h.
\end{align*}
\noindent To bound the interface terms, there is another essential inequality which is very helpful:
\begin{lemma}\cite[Theorem 5.2]{CGPR17} 
$\forall \psi_h\in M_{ph}, \forall\bu_h\in\bV_h $, the following estimate holds true
\begin{align}
|(\psi_h,\bu_h\cdot\bn_{\gamma})_{\Gamma}|\leq C\|\psi_h\|_{M_{1p}}\|\bu_h\|_{L^2(\Omega_f)}.\label{T12}
\end{align}
\end{lemma}
\noindent To analyse the discrete problem (\ref{FNSED1})-(\ref{FNSED2}), we require that the bilinear forms $a_{NS}(\cdot,\cdot)$ and $a_{D}(\cdot,\cdot)$ satisfy adequate continuity and coercivity properties.
\begin{lemma} \label{coercivity} \cite[Lemma $4.8$]{GR09}
There exist constants $C_3,C_4,C_5,C_{6}>0$, independent of $h$ and $\mu$, such that
\begin{align*}
&a_{NS}(\bv_h,\bv_h)\geq C_3\mu\|\bv_h\|_{\bX_1} ^2, \,\, a_{NS}(\bv_h,\bz_h)\leq C_5\|\bv_h\|_{\bX_1} \|\bz_h\|_{\bX_1},\quad \forall \bv_h,\bz_h\in {\bX}_h,\\
&a_D(\psi_{h},\psi_{h})\geq C_{4}\|\psi_{h}\|_{M_{1p}} ^2,\,\,\,\,\, a_D(\phi_{h},\psi_{h})\leq C_6\|\phi_{h}\|_{M_{1p}}\|\psi_{h}\|_{M_{1p}}, \quad \forall \phi_{h},\psi_{h}\in {M_{ph}},
\end{align*}
presuming that the penalty parameter $\sigma_e$ is chosen large enough, when $\epsilon_i$ takes the value $0$ or $-1$.
\end{lemma}
\noindent A uniform discrete inf-sup condition is now presented below.
\begin{lemma}\cite[Proposition $4.13$]{GR09}\label{infsup}
There exists a positive constant $\gamma^{*}$, independent of $h$, such that
\begin{align*}
\inf_{q_h\in M_{fh}} \sup_{\bv_h\in \bV_h} \frac{b_{NS}(\bv_h,q_h)}{\|\bv_h\|_{\bX_1}\|q_h\|_{M_{1f}}}\geq \gamma^{*}.
\end{align*}
\end{lemma}
\noindent Since the interface conditions do not eliminate the nonlinear term from the energy balance, instead of ellipticity, we obtain the following result. 
\begin{lemma} \cite[Lemma $4.11$]{GR09}
The nonlinear term $c_{NS}(\cdot,\cdot;\cdot,\cdot,)$ satisfies 
\begin{align}
c_{NS}(\bz_h,\bz_h;\bv_h,\bv_h)\geq \frac{1}{2}(\bz_h\cdot\bn_{\gamma},\bv_h\cdot\bv_h)_{\Gamma}, \quad \forall \bz_h, \bv_h\in {\bX}_h. \label{nonliearity1}
\end{align}
\end{lemma}
\noindent The next two lemmas help us to handle the nonlinear term efficiently. \begin{lemma}\cite[Lemma $2.3$]{KHR23}
When $\Omega\subset \mathbb{R}^2$ is convex, the following estimate holds true: 
\begin{align}
\|\bw_h\|_{L^4(\Omega)}\leq C \|\bw_h\|^{1/2}_{L^2(\Omega)}\|\bw_h\|^{1/2}_{\bX_1},\quad \forall \bw_h\in\bX_h .\label{nonlinearity8}
\end{align}
\end{lemma}
\begin{lemma} The following inequalities for the trilinear term hold true:
\begin{align}
|&n(\bz_h  ,\bw,\bu_h)|  \leq C\|\bz_h\|_{L^2(\Omega)}|\bw|_{1,4,\Omega_f}\|\bu_h\|_{\bX_1},\,\, \bw\in\textbf{W}^{1,4}(\Omega_f),\,\,\bz_h,\bu_h\in \bV_h,\label{nonlinearity3}\\
|&c_{NS}(\bz,\bz_h;\bu_h,\bv_h)|  \leq C\|\bz_h\|_{\bX_1}\|\bu_h\|_{\bX_1}\|\bv_h\|_{\bX_1}, \,\, \bz\in\bX_1, \bu_h,\bz_h,\bv_h\in\bX_h,\label{nonlinearity4}\\
|&c_{NS}(\bz,\brho  ;\bv,\bu_h)|  \leq C\|\brho\|_{\bX_1}\|\bv\|_{\bX_1}\|\bu_h\|_{\bX_1},\,\, \brho,\bv,\bz\in\bX_1,\,\bu_h\in \bX_h,\label{nonlinearity5}\\
|&l(\brho_h,\bz_h;\bu_h,\bw_h)-l(\bz_h,\bz_h;\bu_h,\bw_h)| \leq C\|\brho_h-\bz_h\|_{\bX_1}\|\bu_h\|_{\bX_1}\|\bw_h\|_{_{L^4(\Omega_f)}}\,,\,\,\brho_h, \bu_h,\bw_h,\bz_h \in \bX_h.\label{nonlinearity7}
\end{align}
\end{lemma}
\begin{proof}
The first two inequalities can be found in \cite[Proposition 4.1, $(4.1)$, $(4.3)$]{GRW005S}. For the third inequality, we refer to \cite[$(3.39)$]{BGR23}.
Also, from \cite[$(3.40)$]{BGR23}, \vspace*{-1mm} 
\begin{align*}
|&l(\brho_h,\bz_h;\bu_h,\bw_h)-l(\bz_h,\bz_h;\bu_h,\bw_h)| \leq C\|\brho_h-\bz_h\|_{L^4(\Omega_f)}\|\bu_h\|_{\bX_1}\|\bw_h\|_{L^4(\Omega_f)}.
\end{align*}
Then using (\ref{T8}), we arrive at \vspace*{-1mm}
\[ |l(\brho_h,\bz_h;\bu_h,\bw_h)-l(\bz_h,\bz_h;\bu_h,\bw_h)| \leq C\|\brho_h-\bz_h\|_{\bX_1}\|\bu_h\|_{\bX_1}\|\bw_h\|_{L^4(\Omega_f)}. \qedhere \] 
\end{proof}
\noindent The existence and uniqueness of the solution of (\ref{FNSED3}) at each time step can be shown as in Section $5.2$ in \cite{CGPR17}. Below in Lemma \ref{apriori}, we present the stability of the scheme. 
\vspace*{-2mm}
\begin{lemma} \label{apriori}
Let $0<\alpha<min\big\{\frac{\mu C_3}{2C_2},\frac{gC_4 }{2C_1}\big\}$. 
If for all $0\leq i\leq N-1$, 
\begin{align}
\|\bU^{i}_h\|_{\bX_1}\leq\frac{\mu C_3}{4C_0}, \label{aprioricondition1}
\end{align}
then the discrete solution $(\{\bU^i_h\}, \{\Phi^i_h\})_{1\geq n\leq N }$ of (\ref{FNSED3}) satisfies
\begin{align*}
\|\bU^n_h\|^2_{L^2(\Omega_f)} +  gS_0\|\Phi^n_h\|^2_{L^2(\Omega_p)} + e^{-2\alpha t_n}\Delta t \sum_{i=1}^{n} e^{2\alpha t_i}\|\bU^i_h\|^2_{\bX_1} + e^{-2\alpha t_n}\Delta t \sum_{i=1}^{n} e^{2\alpha t_i}\|\Phi^i_h\|^2_{M_{1p}}\leq C.
\end{align*}
Here the constant $C=e^{2\alpha \Delta t}\big( \|\bU_h^0\|^2 + \|\Phi_h^0\|^2 \big) + \frac{C_2^2}{\mu C_3}\Delta t \sum_{i=1}^{n} e^{2\alpha t_i} \|\f_f^i\|^2_{L^2(\Omega_f)} + \frac{C_1^2}{2 C_4 g}\Delta t \sum_{i=1}^{n} e^{2\alpha t_i} \|f_p^i\|^2_{L^2(\Omega_p)}. $ 
\end{lemma}
\noindent The above estimate can be obtained by choosing $\bv_h=\bU^i_h$ and $\psi_h=\Phi^i_h$ in (\ref{FNSED3}) and using Lemma \ref{coercivity}.
\vspace*{-3mm}
\section{Approximation Operators}\label{s4.}

In our analysis, we disassociate the pressure from velocity by considering the DG formulation over a weakly divergence-free space. As a result, it is imperative for us to look for an approximation operator from a divergence-free space, say, $J$, to a discrete divergence-free space, where the space $J$ is given by
\[ \textbf{J} = \{\bv\in \bX: \nabla\cdot\bv = 0 ~\text{in}~{\Omega}_f\}. \]
\vspace*{-8mm}
\begin{lemma}\label{iprojection}
There is an approximation $i_h:\bJ\cap \bH^{k_1+1}(\Omega_f) \rightarrow \bV_h$ such that \vspace*{-2mm} 
\begin{align}
\|\bz-i_h\bz\|_{L^2(\Omega_f)} + h\|\bz-i_h\bz\|_{\bX_1} \leq Ch^{k_1+1} |\bz|_{k_1+1,\Omega_f}, \forall\bz\in\bJ\cap\bH^{k_1+1}(\Omega_f).
\end{align}
\end{lemma}
\begin{proof}
Lemma \ref{Piprojection} yields that there exists an operator $\bPi_h:\bX\rightarrow\bX_h$ which satisfies 
\begin{align*}
b_{NS}({\bPi}_h(\bw)-\bw,q_h)=0\,,\quad \forall \bw\in \bX, \,\,\forall q_h\in M_{fh}.
\end{align*}
Restricting ${\bPi}_h$ to $\bJ\cap \bH^{k_1+1}(\Omega_f)$, we observe that $b_{NS}({\bPi}_h(\bw),q_h)=0$, which indicates ${\bPi}_h(\bw)\in\bV_h$. 
This restriction is defined as ${i}_h:\bJ\cap \bH^{k_1+1}(\Omega_f) \rightarrow \bV_h$. From (\ref{Piproperty}), we have 
\begin{equation*}
\|\bz-i_h\bz\|_{L^2(\Omega_f)} + h\,\|\bz-i_h\bz\|_{\bX_1} = \|\bz-\bPi_h\bz\|_{L^2(\Omega_f)} + h\,\|\bz-\bPi_h\bz\|_{\bX_1}\leq Ch^{k_1+1} |\bz|_{k_1+1,\Omega_f}. \hfill \qedhere
\end{equation*}
\end{proof}
\noindent Next, we define an $\bL^2$ projection $\bP_h: \bL^2(\Omega_f)\rightarrow \bV_{h}$: For $\bv\in \bL^2(\Omega_f)$, \vspace*{-1mm}
\begin{align*}
(\bv-\bP_h\bv,{\textbf v}_h)=0,\quad \forall {\textbf v}_h\in \bV_h.
\end{align*}
\noindent 
The following optimal estimate is a consequence of the above lemma, Lemma \ref{iprojection}.
\begin{lemma}\label{pprojection}\cite[Lemma 4.3]{HR82}
The approximation operator $\bP_h$ satisfies \vspace*{-2mm} 
\begin{align}
\|\bw-\bP_h(\bw)\|_{L^2(\Omega_f)}+h\|\bw-\bP_h(\bw)\|_{\bX_1}\leq Ch^{k_1+1} |\bw|_{k_1+1,\Omega_f}\quad\forall \bw\in\bJ\cap\bH^{k_1+1}(\Omega_f) .\label{projection2}
\end{align}
\noindent presuming the constant $C>0$ is independent of $h$.
\end{lemma}
\noindent The stability of the projection $\bP_h$ can be obtained using the triangle inequality, the Lagrange interpolant, (\ref{T4}), and the projection property (\ref{projection2}): \vspace*{-2mm}
\begin{align}
\|\bP_h\bv\|_{\bX_f}\leq C|\bv|_{1,\Omega_f},\quad \forall\,\bv\in \bJ.\label{Ph}
\end{align} 
\noindent We observe that higher-order polynomial approximation requires sufficient regularity for the weak solutions 
$(\bu, p,\phi)$, and therefore, we assume the following:\\ 
\noindent $(\textbf{A}\textbf{1}).$ The exact solutions satisfy:\\ \begin{tabular}{c c c }
	&$\bu \in L^{\infty}(0\,,T;\bH^{k_1+1}(\Omega_f))\cap L^{2}(0\,,T;\bH^{k_1+1}(\Omega_f) )$, $p \in L^{\infty}(0\,,T;\bH^{k_1}(\Omega_f))\cap L^{2}(0\,,T;\bH^{k_1}(\Omega_f) )$     \\
	&$\phi  \in L^{\infty}(0\,,T;H^{k_2+1}(\Omega_p))\cap L^{2}(0\,,T;H^{k_2+1}(\Omega_p) )$,\,\, $p_t\in L^2(0,T;H^{k_1}(\Omega_f))$ ,        \\
   & $\bu_t \in L^{2}(0\,,T;\bH^{k_1+1}(\Omega_f) )\cap L^{\infty}(0\,,T;\bH^{k_1}(\Omega_f))$,\,\, $\bu_{tt} \in L^{2}(0\,,T;\bL^2(\Omega_f) )$,  \\
 &  $\phi_t \in L^{2}(0\,,T;H^{k_2+1}(\Omega_p) )\cap L^{\infty}(0\,,T;H^{k_2}(\Omega_p))$,\,\, $\phi_{tt}  \in L^{2}(0\,,T;L^2(\Omega_p) )$. \\
		\end{tabular} \vspace*{3mm}
		
\noindent
We now present below the modified Oseen operators, which are crucial for our $L^2$-error analysis. This and the rest of our analysis will be carried out for the two-dimensional SIPG case, as NIPG and IIPG cases will lead us to suboptimal results; for the three-dimensional SIPG case, please see Remark \ref{3D}.
\subsection{Modified Oseen Operators}\label{s5.} 
As mentioned earlier, the NSD problem offers a unique difficulty due to the presence of the interface terms, and the standard operators have not been adequate to handle them in the case of $L^\infty(\textbf{L}^2)$ estimates. Therefore, we present below the modified Oseen type operators $\bO_{fh}$ for velocity and $O_{ph}$ for piezometric head, respectively, satisfying \vspace*{-2mm}
 \begin{align}
a_{NS}(\bu-\bO_{fh}\bu,\bv_h)+ ga_D(\phi-O_{ph} \phi,\psi_h)+  g (\phi-O_{ph} \phi,\bv_h\cdot{\bn_{\gamma}})_{\Gamma} - g((\bu-\bO_{fh}\bu)\cdot\bn_{\gamma},\psi_h)_{\Gamma} \nonumber \\
+ \frac{1}{G}((\bu -\bO_{fh}\bu)\cdot\btau_{\gamma},\bv_h\cdot\btau_{_{\gamma}})_{\Gamma} + c_{NS}(\bu  ,\bu;\bu-\bO_{fh}\bu,\bv_h)+ c_{NS}(\bu  ,\bu-\bO_{fh}\bu;\bu,\bv_h)\nonumber\\
 =\,- b_{NS}(\bv_h,p), \forall \bv_h\in\bV_h,\forall \psi_h\in M_{ph}, \label{Oseenprojection}
\end{align}
given ($\bu$, $p$, $\phi$) is the weak solution of the NSD problem. \\
Being a linearized steady NSD problem, the existence and uniqueness of (\ref{Oseenprojection}) can be proved in similar ways as in section $5$ of \cite{GR09}.
Approximation properties of the operators $\bO_{fh}$ and $O_{ph}$ follow next.
\begin{lemma} \label{Opproperty}
If $\bu$ satisfies
\begin{align}
\|\textbf{D}(\bu)\|_{L^2(\Omega_f)} \leq \frac{C_3\mu}{4(1+d_1d_2)}, \label{os}
\end{align}
then terms $\bu-\bO_{fh}\bu$ and $\phi-O_{ph}\phi$ satisfy the following estimates:
\begin{align}
\|\bu-\bO_{fh}\bu\|^2_{L^2(\Omega_f)}+h^2\|\bu&-\bO_{fh}\bu\|^2_{\bX_1}+\|\phi-O_{ph}\phi\|^2_{L^2(\Omega_p)} +h^2\|\phi-O_{ph}\phi\|^2_{M_{1p}}\nonumber\\
&\leq C (h^{2k_1+2}(|\bu|^2_{k_1+1,\Omega_f}
+|p|^2_{k_1,\Omega_f})+h^{2k_2+2}|\phi|^2_{k_2+1,\Omega_p}),\label{Oseenprojection1}\\
\|(\bu-\bO_{fh}\bu)_t\|^2_{L^2(\Omega_f)}+h^2\|(\bu&-\bO_{fh}\bu)_t\|^2_{\bX_1}+\|(\phi-O_{ph}\phi)_t\|^2_{L^2(\Omega_p)} +h^2\|(\phi-O_{ph}\phi)_t\|^2_{M_{1p}}\nonumber\\
&\leq C (h^{2k_1+2}(|\bu_t|^2_{k_1+1,\Omega_f}
+|p_t|^2_{k_1,\Omega_f})+h^{2k_2+2}|\phi_t|^2_{k_2+1,\Omega_p}),\label{Oseenprojection2}
\end{align}
where $C = C(\mu,g, C_3, C_4, {d_1}, d_2 )>0$ is independent of $h$.
\end{lemma}
\begin{proof}
Since
\vspace{-2mm}
 \begin{align}
\|\bu-\bO_{fh}\bu\|_{\bX_1} &\leq \|\bu-\bP_{h}\bu\|_{\bX_1}+ \|\bP_h \bu- \bO_{fh}\bu\|_{\bX_1}, \label{O1}\\
\|\phi-O_{ph}\phi\|_{M_{1p}} &\leq \|\phi-L_h \phi\|_{M_{1p}}+\|L_h \phi-O_{ph}\phi\|_{M_{1p}}, \label{O2}
\end{align}
it is sufficient to estimate $\bP_h \bu- \bO_{fh}\bu $ and $L_h \phi-O_{ph}\phi $. We choose $ \bv_h=\bP_h\bu-\bO_{fh}\bu=:\U$ and $\psi_h= L_h\phi-O_{ph}\phi=:\mathcal{J}$ in (\ref{Oseenprojection}) and observe that, since $b_{NS}(\U,r_h)=0$,
\begin{align}
 a_{NS}(\U,\U)+ \frac{1}{G}(\U \cdot\btau_{\gamma},\U\cdot\btau_{\gamma})_{\Gamma}
+ ga_D(\mathcal{J} ,\mathcal{J}  )  
\leq    - b_{NS}&(\U,r_h p-p) - a_{NS}(\bu-\bP_h \bu,\U)  
- ga_D(\phi-L_h \phi ,\mathcal{J} ) \nonumber \\
 - \frac{1}{G}((\bu-\bP_h \bu) \cdot\btau^j_{\gamma},\U\cdot\btau^j_{\gamma})_{\Gamma}  -  g ( \phi-L_h\phi,\U \cdot{\bn_{\gamma}})_{\Gamma}
 &+  g((\bu-\bP_{h}\bu)\cdot\bn_{\gamma},\mathcal{J} )_{\Gamma}
  - c_{NS}(\bu  ,\bu; \bu-\bP_h \bu,\U)\nonumber \\ -   c_{NS}(\bu  ,\bu;\U,\U)
 - c_{NS}(\bu  ,\bu -\bP_h \bu;\bu,\U)
  - c_{NS}(&\bu  ,\U;\bu,\U)\quad = \quad \sum_{i=1}^{10} A_i 
\label{O3}.
\end{align}
\noindent Now using Lemma \ref{coercivity}, the left side of (\ref{O3}) can be bounded below by
\begin{align}
 \mu C_3 \|\U\|^2_{\bX_1} + C_4 g \|\mathcal{J}\|^2_{M_{1p}} +\frac{1}{G}\|\U\cdot\btau_{\gamma}\|^2_{L^2(\Gamma)}.\label{O4}    
\end{align}%
For the right hand side of (\ref{O3}), we utilize the standard Cauchy-Schwarz and Young's inequalities, the projection properties (\ref{projection1}), (\ref{projection3}), (\ref{projection2}), the trace inequalities (\ref{T1}) and (\ref{T4}), and obtain
\begin{align*}
&\big|A_1\big| = \Big| -\sum_{E\in\mathcal{E}_{fh}} (r_h p-p,\nabla\cdot\U)_E + \sum_{e\in \Gamma_{fh}\cup\Gamma_f} (\{r_h p-p\},[\U]\cdot \bn_e)_e \Big| 
\leq \frac{C_3\mu}{48}\|\U\|^2_{\bX_1} + Ch^{2k_1}|p|^2_{k_1,\Omega_f}.\\ 
&\big|A_2\big|+\big|A_3\big|\leq \frac{C_3\mu}{48}\|\U\|^2_{\bX_1}  + \frac{gC_4}{10}\|\mathcal{J}\|^2_{\bM_{1p}}+ C\big(h^{2k_1}|\bu|^2_{k_1+1,\Omega_f}+ h^{2k_2}|\phi|^2_{k_2+1,\Omega_p}\big).
\end{align*}%
Additionally, we apply (\ref{T9}) to bound the interface terms as 
\begin{align*} 
&\big|A_4\big|\leq \frac{1}{G} \|\bu-\bP_h \bu\|_{L^2(\Gamma)}\|\U\cdot\btau_{\gamma}\|_{L^2(\Gamma)}\nonumber\\
&\qquad\leq \frac{C}{G} \bigg(  \sum_{E\in\mathcal{E}_{fh}} h_E^{-1}\big(\|\bu-\bP_h \bu\|^2_{L^2(E)} + h_E^2 \|\nabla(\bu-\bP_h \bu)\|^2_{L^2(E)}\big) \bigg)^{1/2} \|\U\cdot\btau_\gamma\|_{L^2(\Gamma)}\nonumber\\
&\qquad\leq \frac{1}{2G}\|\U\cdot\btau_\gamma\|^2_{L^2(\Gamma)} + Ch^{2k_1+1}|\bu|^2_{k_1+1,\Omega_f},\\
&\big|A_5\big|+\big|A_6\big| \leq \frac{C_3\mu}{48}\|\U\|^2_{\bX_1}+ \frac{gC_4}{10}\|\mathcal{J}\|^2_{M_{1p}} + C\big(h^{2k_2+1}|\phi|^2_{k_2+1,\Omega_p}  + h^{2k_1+1}|\bu|^2_{k_1+1,\Omega_f}\big).
\end{align*}%
Next, we bound the nonlinear terms using (\ref{T9}), (\ref{nonliearity1}), and (\ref{os}) to arrive at 
\begin{align*}
\big| A_8\big| & =\Big|\sum_{E\in\mathcal{E}_{fh}}  (\bu\cdot \nabla \U,\U)_E +\frac{1}{2}\sum_{E\in \mathcal{E}_{1h}}(\nabla\cdot \bu,\U\cdot \U)  -\frac{1}{2}\sum_{e\in \Gamma_{fh}\cup\Gamma_f}([\bu]\cdot \bn_e\{\U \cdot \U\}) \\
  & \hspace*{1cm}+  \sum_{E\in\mathcal{E}_{fh}}(|\{\bu\}\cdot \bn_E|[\U], \U^{int})_{\partial E_-(\bu) \setminus \Gamma} \Big| \\
  & \leq \|\bu\|_{L^4(\Omega_f)}\|\U\|^2_{\bX_1}\leq(1+{d_2}d_1)\|\textbf{D}(\bu)\|_{L^2(\Omega_f)}\|\U\|^2_{\bX_1}\leq \frac{C_3\mu}{4}\|\U\|^2_{\bX_1},\\
\big|A_{10}\big| & \leq \frac{C_3\mu}{4}\|\U\|^2_{\bX_1}.
\end{align*}
\noindent Also, the rest of the nonlinear terms can be bounded using (\ref{nonlinearity5}) along with Young's inequality as
\begin{align*}
\big|A_7\big| + \big| A_9\big|&\leq C \|\bu\|_{2,\Omega_f} \|\bu-\bP_h\bu\|_{\bX_1}\|\U\|_{\bX_1}\leq \frac{C_3\mu}{48}\|\U\|^2_{\bX_h} + Ch^{2k_1}|\bu|^2_{k_1+1,\Omega_f}.
\end{align*}%
Incorporating (\ref{O4}) and the above bounds in (\ref{O3}), we obtain
\begin{align*}
\mu C_3 \|\U\|^2_{\bX_1} + gC_4 \|\mathcal{J}\|^2_{M_{1p}} +\frac{1}{G}\|\U\cdot\btau_\gamma\|^2_{L^2(\Gamma)}
 \leq  C(h^{2k_1}(|\bu|^2_{k_1+1,\Omega_f} 
+|p|^2_{k_1,\Omega_f } )+  h^{2k_2}|\phi|^2_{k_2+1,\Omega_p}).
\end{align*}%
\noindent Since $\bu-\bO_{fh}\bu=\bu-\bP_h\bu+\U$ and $\phi-O_{ph}\phi=\phi-L_h\phi+\mathcal{J}$, using approximation properties (\ref{projection3}) and (\ref{projection2}), we complete the energy norm estimates of $\bu-\bO_{fh}\bu$ and $\phi-O_{ph}\phi$:
\begin{align}
\|\bu-\bO_{fh}\bu\|^2_{\bX_1}+\|\phi-O_{ph}\phi\|^2_{M_{1p}}
\leq C (h^{2k_1}(|\bu|^2_{k_1+1,\Omega_f}
+|p|^2_{k_1,\Omega_f})+h^{2k_2}|\phi|^2_{k_2+1,\Omega_p}).\label{O18}
\end{align}%
\noindent For the $L^2$-norm estimate, we apply the Aubin-Nitsche duality argument. We consider the following dual problem: Given $\textbf{g}_1\in L^2(\Omega_f)^2$, $g_2\in L^2(\Omega_p)$, find $(\bz,\varsigma,\nu)\in \bJ\times M_f\times M_p $ such that $\forall\,\,(\bv,q,\psi)\in \bX\times M_f\times M_p$,
\begin{align}
 2\mu(\bD(\bv),\bD(\bz)) + g\bK (\nabla\theta,\nabla\nu) + (\bu\cdot\nabla\bv,\bz) +  (\bv\cdot\nabla\bu,\bz) + \frac{1}{G}(\bv\cdot\btau_{\gamma},\bz\cdot\btau_{\gamma})_{\Gamma}\nonumber\\
 +   g(\bv\cdot\bn_{\gamma},\nu)_{\Gamma} -  g(\psi,\bz\cdot\bn_{\gamma})_{\Gamma} + (\varsigma, \nabla\cdot\bv) = (\textbf{g}_1,\bv)_{\Omega_f}+(g_2,\psi)_{\Omega_p}\label{Dual}
 \end{align}
\noindent Following the argument in \cite[pp. 3332]{CMX09}, one finds that the two trilinear terms in the dual problem (\ref{Dual}) are bounded for a sufficiently large $\mu$, and the regularity assumption in $(\textbf{A1})$. Consequently, we may assume that the solution of the dual problem satisfies the following regularity result:
\begin{align}
\bigg(\|\bz\|^2_{2,\Omega_f}+\|\nu\|^2_{2,\Omega_p}\bigg)^{1/2}+ \|\varsigma\|_{1,\Omega_f}\leq \|\textbf{g}_1\|_{L^2(\Omega_f)}+\|g_2\|_{L^2(\Omega_p)}.\label{Dualregularity}
\end{align}
\noindent Due to consistency, one can rewrite (\ref{Dual}) as follows: \vspace{-1mm}
\begin{align*}
 a_{NS}(\bv,\bz) &+ g a_D(\psi,\nu) + c_{NS}(\bu,\bu;\bv,\bz) +  c_{NS}(\bu,\bv;\bu,\bz) + \frac{1}{G}(\bv\cdot\btau_{\gamma},\bz\cdot\btau_{\gamma})_{\Gamma}
 +   g(\bv\cdot\bn_{\gamma},\nu)_{\Gamma}\nonumber\\
  &-  g(\psi,\bz\cdot\bn_{\gamma})_{\Gamma} + b_{NS}(\bv,\varsigma) = (\textbf{g}_1,\bv)_{\Omega_f}+(g_2,\psi)_{\Omega_p},\quad \forall\bv\in\bX_1, \psi\in M_{1p}.
\end{align*}
\noindent Allowing $\textbf{g}_1=\bu-\bO_{fh}\bu$, $g_2=O_{ph}\phi-\phi$, we put $\bv=\bu-\bO_{fh}\bu $ and $\psi = O_{ph }\phi-\phi$ in the above equation. We then subtract the resulting equation from (\ref{Oseenprojection}) after setting $\bv_h= \bP_h \bz$ and $\psi_h=-L_h \nu$ to obtain
\begin{align}
&\|\bu-\bO_{fh}\bu\|^2_{L^2(\Omega_f)}+\|O_{ph }\phi-\phi\|^2_{L^p(\Omega_p)} = a_{NS}(\bu-\bO_{fh}\bu,\bz-\bP_h \bz) + ga_D(\phi-O_{ph }\phi,L_h \nu-\nu)\nonumber\\ 
 &\quad\quad+ b_{NS}(\bu-\bO_{fh}\bu, \varsigma)+b_{NS}(\bz-\bP_h \bz, p-r_h p) +  g ((\bu-\bO_{fh}\bu)\cdot{\bn_{\gamma}},\nu-L_h \nu)_{\Gamma}\nonumber\\
& \quad\quad - g (\phi-O_{ph }\phi,(\bz-\bP_h \bz) \cdot{\bn_{\gamma}})_{\Gamma}+\frac{1}{G}((\bu-\bO_{fh}\bu)\cdot{\btau_{\gamma}},(\bz-\bP_h \bz)\cdot{\btau_{\gamma}})_{\Gamma}\nonumber\\
&\quad\quad +c_{NS}(\bu,\bu ;\bu-\bO_{fh}\bu,\bz-\bP_h \bz) +c_{NS}(\bu,\bu-\bO_{fh}\bu;\bu,\bz-\bP_h \bz)=\sum_{i=1}^{9}B_i.\label{O21}
\end{align}
\noindent Here we have used the fact $[\bz]\cdot\bn_e=0$ in the fourth term on the right side of (\ref{O21}).
\noindent The first two bilinear terms can be bounded using similar techniques applied for the second and third terms in (\ref{O3}) as follows:  {\small %
\begin{align*}
|B_1+B_2|
\leq\frac{C_3}{36} \|\bu-\bO_{fh}\bu\|^2_{L^2(\Omega_f)} + \frac{C_4}{36} \|\theta-L_h \theta\|^2_{L^2(\Omega_p)}+ C\big(h^{2k_1+2}(|\bu|^2_{k_1+1,\Omega_f}+ |p|^2_{k_1,\Omega_f }) 
 +  h^{2k_2+2}|\phi|^2_{k_2+1,\Omega_p}\big).
\end{align*}}%
\noindent Using the definition of $\bV_h$, one finds
\begin{align*}
|B_3|&= \mid b_{NS}(\bu-\bO_{fh}\bu, \varsigma-r_h\varsigma) \mid = \mid b_{NS}(\bu-\bP_h\bu, \varsigma-r_h\varsigma)+b_{NS}(\bP_h\bu-\bO_{fh}\bu, \varsigma-r_h\varsigma) \mid.
\end{align*} 
\noindent Now expanding the terms and using standard inequalities with the projection property (\ref{projection1}), we find 
\begin{align*}
&|B_3|\leq \Big| -\sum_{E\in\mathcal{E}_{fh}} (\varsigma-r_h\varsigma,\nabla\cdot (\bu-\bP_h\bu))_E + \sum_{e\in \Gamma_{fh}\cup\Gamma_f} (\{\varsigma-r_h\varsigma\},[\bu-\bP_h\bu]\cdot \bn_e)_e \Big| \nonumber\\
&\qquad+ 
  \Big| \sum_{e\in \Gamma_{fh}\cup\Gamma_f} (\{\varsigma-r_h\varsigma\},[\bP_h\bu-\bO_{fh}\bu]\cdot \bn_e)_e\Big|\nonumber\\
&\leq\frac{C_3\mu}{36} \|\bu-\bO_{fh}\bu\|^2_{L^2(\Omega_f)} + \frac{C_4}{36} \|\theta-L_h \theta\|^2_{L^2(\Omega_p)} +  C\big(h^{2k_1+2}(|\bu|^2_{k_1+1,\Omega_f}+ |p|^2_{k_1,\Omega_f }) 
 +  h^{2k_2+2}|\phi|^2_{k_2+1,\Omega_p}\big).
\end{align*}
Similarly, the fourth term can be bounded as
\begin{align*}
|B_4| &\leq\frac{C_3\mu}{36} \|\bu-\bO_{fh}\bu\|^2_{L^2(\Omega_f)} + \frac{C_4}{36} \|\theta-L_h \theta\|^2_{L^2(\Omega_p)} + C\big(h^{2k_1+2} |p|^2_{k_1,\Omega_f }).
\end{align*}%
\noindent For the fifth term, applying (\ref{Poincare}), (\ref{projection3}), (\ref{T1}), and (\ref{O18}), we obtain
\begin{align*}
|B_5|&\leq  \|(\bu-\bO_{fh}\bu)\cdot{\bn}\|_{L^2(\Gamma)}\|L_h \nu-\nu\|_{L^2(\Gamma)}\nonumber\\
&\leq C \bigg(  \sum_{E\in\mathcal{E}_{fh}} h_E^{-1}\big(\|\bu-\bO_{fh}\bu\|^2_{L^2(E)} + h_E^2 \|\nabla(\bu-\bO_{fh}\bu)\|^2_{L^2(E)}\big) \bigg)^{1/2}\\
&\qquad \times \bigg(  \sum_{E\in\mathcal{E}_{fh}} h_E^{-1}\big(\|L_h \nu-\nu\|^2_{L^2(E)} + h_E^2 \|\nabla(L_h \nu-\nu)\|^2_{L^2(E)}\big) \bigg)^{1/2}\nonumber\\
&\leq  \frac{C_3\mu}{36} \|\bu-\bO_{fh}\bu\|^2_{L^2(\Omega_f)} + \frac{C_4}{36} \|\theta-L_h \theta\|^2_{L^2(\Omega_p)} + Ch^{2k_1+2}|\bu|^2_{k_1+1,\Omega_f}+ Ch^{2k_2+2}|\phi|^2_{k_2+1,\Omega_p}.
\end{align*}
\noindent Similar to the above term, the remaining interface terms in (\ref{O21}) can be bounded as follows:
\begin{align*}
|\sum_{i=6}^7B_i|
  \leq \frac{C_3\mu}{36} \|\bu-\bO_{fh}\bu\|^2_{L^2(\Omega_f)} + \frac{C_4}{36} \|\theta-L_h \theta\|^2_{L^2(\Omega_p)} + Ch^{2k_1+2}|\bu|^2_{k_1+1,\Omega_f}+ Ch^{2k_2+2}|\phi|^2_{k_2+1,\Omega_p}.
\end{align*}
\noindent The nonlinear terms are bounded using the Cauchy-Schwarz inequality, (\ref{O18}), and Young's inequality as
\begin{align*}
& |B_8+B_9|
\leq \|\bu\|_{2,\Omega_f}\|\bu-\bO_{fh}\bu\|_{\bX_h}\|\bw-\bP_h \bw \|_{\bX_1}\nonumber\\
&\,\,\leq\frac{C_3\mu}{36} \|\bu-\bO_{fh}\bu\|^2_{L^2(\Omega_f)} + \frac{C_4}{36} \|\theta-L_h \theta\|^2_{L^2(\Omega_p)}+C (h^{2k_1+2}(|\bu|^2_{k_1+1,\Omega_f}
+|p|^2_{k_1,\Omega_f})+h^{2k_2+2}|\phi|^2_{k_2+1,\Omega_p}).
\end{align*}
\noindent Inserting all the above bounds in (\ref{O21}), we achieve the desired result.\par
\noindent Repeating the aforementioned steps, we obtain the estimates (\ref{Oseenprojection2}) consisting of $(\bu-\bO_{fh}\bu)_t$ and $(\phi-O_{ph }\phi)_t$. Here, in place of equation (\ref{O3}), we utilize the differentiated one, in time. Substitute $\bv_h =(\bP_h \bu- \bO_{fh}\bu)_t$ and $\psi_h=(L_h \phi-O_{ph}\phi)_t $ in the resulting equation, and finally, we choose the right side as $(\bu-\bO_{fh}\bu)_t$ and $(L_h \phi-O_{ph}\phi)_t$ in the dual problem. This concludes the proof.
\end{proof}

\section{Fully Discrete Error Analysis}\label{s6.}
This section deals with the error bounds of the fully discrete scheme. 
The NS velocity and piezometric head errors are denoted by $\e_f^i=\bu^i-\bU_h^i$ and $e_p^i=\phi^i-\Phi_h^i$ for $i=1,2,...,N$, respectively.
Due to the consistency of the numerical scheme (\ref{FNSED1})-(\ref{FNSED2}), we observe that the weak solution $(\bu,p,\phi)$ satisfies (\ref{FNSED1}) with $\partial_t\bU_h^i$ replaced by $\bu_t^i$ at $t=t_i$. 
Therefore, subtracting (\ref{FNSED3}) from the resulting equation, we obtain the error equation such that $\forall(\bv_h,\psi_h)\in(\bV_h,M_{ph})$,
\begin{align}
&(\partial_t \e_f^i,\bv_h)+ gS_0(\partial_t e_p^i,\psi_h)+a_{NS}(\e_f^i,\bv_h)+ ga_D(e_p^i,\psi_h) +\frac{1}{G}(\e_f^i\cdot\btau_{\gamma},\bv_h\cdot\btau_{\gamma})_{{\Gamma}}\nonumber\\
 &\quad +c_{NS}(\bu^i,  \bu^i;\bu^i,\bv_h)
-c_{NS}(\bU_h^{i-1},  \bU_h^{i-1};\bU_h^{i},\bv_h )
+  g(e_p^i,\bv_h\cdot{\bn_{\gamma}})_{{\Gamma}}  - g(\e_f^i\cdot\bn_{\gamma},\psi_h)_{{\Gamma}}  \nonumber\\
&\qquad +  b_{NS}(\bv_h,p^i)=(\partial_t \bu^i-\bu_t^i,\bv_h)  +  gS_0(\partial_t \phi^i-\phi_t^i,\psi_h).\label{EFNSED9}
\end{align}
We observe that the error due to the first-order approximation of the derivative for any function $\nu(t)$ is
\begin{align}
(\partial_t \nu^i,v)- ( \nu_t^i,v) = -\frac{1}{\Delta t} \int_{t_i-1}^{t_i} (s-t_{i-1})(\nu_{ss}(s),v)ds, \label{Taylor}
\end{align}
where $\nu_t^i=\nu_t(x_i)$.
For our subsequent analysis, we split the errors as \vspace*{-3mm}
\begin{align*}
\e_f^i&=\bu^i-\bO_{fh}\bu^i+\bO_{fh}\bu^i-\bU_h^i=\bzeta^i+\bchi^i,\\ 
e_p^i&=\phi^i-O_{ph}\phi^i+O_{ph}\phi^i-\Phi_h^i=\eta^i+\xi^i.
\end{align*}%
By means of the split in the errors, we rewrite the error equation (\ref{EFNSED9}) as
\begin{align}
&(\partial_t \bchi^i,\bv_h)+ gS_0(\partial_t \xi^i,\psi_h)+a_{NS}(\bchi^i,\bv_h)+ ga_D(\xi^i,\psi_h) +\frac{1}{G}(\bchi^i\cdot\btau_{\gamma},\bv_h\cdot\btau_{\gamma})_{{\Gamma}}+  g(\xi^i,\bv_h\cdot{\bn_{\gamma}})_{{\Gamma}}\nonumber\\
& \quad- g(\bchi^i\cdot\bn_{\gamma},\psi_h)_{{\Gamma}}+c_{NS}(\bu^i,  \bu^i;\bu^i,\bv_h)
- c_{NS}(\bU_h^{i-1},  \bU_h^{i-1};\bU_h^{i},\bv_h )
  =(\partial_t \bu^i-\bu_t^i,\bv_h)\nonumber\\
 &\quad-(\partial_t \bzeta^i,\bv_h)  +  gS_0(\partial_t \phi^i-\phi_t^i,\psi_h)    - gS_0(\partial_t \eta^i,\psi_h)-  a_{NS}(\bzeta^i,\bv_h) - ga_D(\eta^i,\psi_h)  -  g( \eta^i,\bv_h\cdot{\bn_{\gamma}})_{{\Gamma}}\nonumber\\
 &\quad+   g(\bzeta^i\cdot\bn_{\gamma},\psi_h)_{{\Gamma}}- \frac{1}{G}(\bzeta^i\cdot\btau_{\gamma},\bv_h\cdot\btau_{\gamma})_{{\Gamma}} -  b_{NS}(\bv_h,p^i) \quad \forall \bv_h\in\bV_h,\psi_h\in M_{ph}.\label{EFNSED10}
\end{align}
\subsection{Error Analysis for Velocity and Piezometric Head } \label{s6}
We present here the optimal estimates for $\bchi^i$ and $\xi^i$ in the $L^{\infty}(\textbf{L}^2)$ and $L^{\infty}(L^2)$-norms, respectively. Before proceeding, we take note of the fact that the constants appearing in the subsequent results depend exponentially on time, although independent of $h$ and $E$, and hence will be denoted by $K_T$.
\begin{lemma}\label{chixi}
Let the assumption $(\textbf{A}\textbf{1})$ and the hypothesis of the Lemma \ref{apriori} be satisfied. Furthermore, assume that the initial values of discrete velocity and piezometric head are $\bU_h^0=\bP_h\bu_0$ and $\Phi_h^0=L_h\phi_0$, respectively. Then for $K_T > 0$, the errors $\bchi^i=\bO_{fh}\bu^i-\bU_h^i$ and $\xi^i=O_{ph}\phi^i-\Phi_h^i$ satisfy
 \begin{align}
\|\bchi^n\|^2_{L^2(\Omega_f)} +  gS_0\|\xi^n\|^2_{L^2(\Omega_p)} + \mu C_3 e^{-2\alpha t_n}&\Delta t \sum_{i=1}^{n} e^{2\alpha t_n}\|\bchi^i\|^2_{\bX_1} +  gC_4 e^{-2\alpha t_n}\Delta t \sum_{i=1}^{n} e^{2\alpha t_n}\|\xi^i\|^2_{M_{1p}}\nonumber\\
+ e^{-2\alpha t_n}\Delta t  \sum_{i=1}^{n} e^{2\alpha t_n}\frac{1}{{G}}\|\bchi^i\cdot\btau_{\gamma}\|^2_{L^2(\Gamma)} &\leq K_T(h^{2k_1+2}+h^{2k_2+2}+\Delta t^2). \label{erroresult1}
\end{align}
\end{lemma}%
\begin{proof}
We now use the definition of modified Oseen operators (\ref{Oseenprojection}) in (\ref{EFNSED10}) and put $\bv_h=\bchi^i$, $\psi_h=\xi^i$. Since modified Oseen operators absorb the interface terms, they will not appear on the right hand side of the resulting equation. Next, we observe that\vspace*{-2mm}
\begin{align*}
(\partial_t \bchi^i,\bchi^i)\geq \frac{1}{2} \partial_t\| \bchi^i\|^2_{L^2(\Omega_f)}, \quad  (\partial_t \xi^i,\xi^i)\geq \frac{1}{2} \partial_t\| \xi^i\|^2_{L^2(\Omega_p)}.
\end{align*}
Next, an application of (\ref{Taylor}), the Cauchy-Schwarz and Young's inequalities leads to 
\begin{align*}
\Big| (\partial_t \bu^i - \bu_t^i,\bchi^i) \Big| &
\leq \frac{\mu C_3}{40} \|\bchi^i\|^2_{\bX_1} + C\Delta t  \int_{t_{i-1}}^{t_i} \|\bu_{tt}\|^2_{L^2(\Omega_f)}  ds,\\ 
\mid  (\partial_t \phi^i-\phi_t^i,\xi^i) \mid 
&\leq \frac{gC_4}{4} \|\xi^i\|^2_{M_{1p}} + C\Delta t  \int_{t_{i-1}}^{t_i} \|\phi_{tt}\|^2_{L^2(\Omega_p)}  ds. 
\end{align*}
And since
$
 \partial_t \bzeta^i   =  \frac{1}{\Delta s} \int_{t_{i-1}}^{t_i} \frac{\partial}{\partial t} (\bzeta^i(s))ds,
$
based on (\ref{Oseenprojection1}) and the standard inequalities, we find that
\begin{align*}
\mid (\partial_t \bzeta^i,\bchi^i)  \mid &= \Big| \frac{1}{\Delta t} \int_{t_{i-1}}^{t_i} (\frac{\partial}{\partial t} (\bzeta^i(s)),\bchi^i)ds\Big| \,\,
 \leq \frac{\mu C_3}{40}\|\bchi^i\|^2_{\bX_1} +  C\frac{h^{2k_1+2}}{\Delta t} |\bu_t|^2_{k_1+1,\Omega_f},\\
 \mid (\partial_t \eta^i,\xi^i) \mid &= \Big| \frac{1}{\Delta t} \int_{t_{i-1}}^{t_i} (\frac{\partial}{\partial t} (\eta^i(s)),\xi^i)ds\Big| \,\,
 \leq \frac{gC_4}{4}\|\xi^i\|^2_{M_{1p}} +  C\frac{h^{2k_2+2}}{\Delta t} |\phi_t|^2_{k_1+1,\Omega_p}.
\end{align*}
\noindent
Incorporating all the above bounds along with (\ref{EFNSED13}), multiplying by $\Delta t e^{2\alpha t_i}$, and finally summing the resulting inequality from $i=1\,\,\text{to}\,\,n ~(\leq N)$ we arrive at
\begin{align}
& e^{2\alpha t_n} \|\bchi^n\|^2_{L^2(\Omega_f)} +  gS_0 e^{2\alpha t_n} \|\xi^n\|^2_{L^2(\Omega_p)}  +  \mu C_3  \Delta t \sum_{i=1}^{n}  e^{2\alpha t_i} \|\bchi^{i}\|^2_{\bX_1}  
+ gC_4\Delta t \sum_{i=1}^{n}  e^{2\alpha t_i} \|\xi^i\|^2_{M_{1p}} \nonumber\\
&\,\,+  \Delta t  \sum_{i=1}^{n} e^{2\alpha t_i}\frac{1}{{G}}\|\bchi^i\cdot\btau_{\gamma}\|^2_{L^2(\Gamma)}
\leq C {\Delta t}\big( \sum_{i=1}^{n}  e^{2\alpha t_i} \|\bchi^i\|^2_{L^2(\Omega_f)} +  gS_0 \sum_{i=1}^{n}  e^{2\alpha t_i}  \|\xi^i\|^2_{L^2(\Omega_p)} \big )\nonumber\\
&\quad\quad\quad\quad\quad\quad\quad\quad\quad+ \Delta t \sum_{i=1}^{n}  e^{2\alpha t_i}\bigg(\sum_{i=1}^{9} H_i\bigg) +  C(h^{2k_1+2}+h^{2k_2+2})e^{2\alpha t_n} + 
C {\Delta t}^2 e^{2\alpha t_{n+1}}\,. \label{EFNSED26}
\end{align}\vspace*{-2mm}
Note that the $H_i$'s appearing above are the trilinear terms, once we rewrite them as follows:
\begin{align}
&\ c_{NS}(\bU_h^{i-1},  \bU_h^{i-1};\bU_h^{i},\bchi^i )  - c_{NS}(\bu^i,  \bu^i;\bu^i,\bchi^i) + c_{NS}(\bu^{i}, \bu^i; \bzeta^{i},\bchi^i )  + c_{NS}(\bu^{i}, \bzeta^{i}; \bu^i ,\bchi^i )\nonumber\\
& = - c_{NS}(\bU_h^{i-1},  \bU_h^{i-1};\bchi^i,\bchi^i )
  - c_{NS}(\bU_h^{i-1},  \bchi^{i-1};\bu^i,\bchi^i )
  + c_{NS}(\bU_h^{i-1},  \bchi^{i-1};\bzeta^i,\bchi^i )\nonumber\\
  &\quad + c_{NS}(\bU_h^{i-1},  \bzeta^{i-1};\bzeta^i,\bchi^i )+ c_{NS}(\bU_h^{i-1}, \bu^i- \bu^{i-1};\bu^i,\bchi^i )
  - c_{NS}(\bU_h^{i-1}, \bzeta^i- \bzeta^{i-1};\bu^i,\bchi^i )\nonumber\\
   &\quad + c_{NS}(\bU_h^{i-1}, \bu^i- \bu^{i-1};\bzeta^i,\bchi^i )+ c_{NS}(\bu^{i}, \bu^i; \bzeta^{i},\bchi^i )
  - c_{NS}(\bU_h^{i-1}, \bu^i; \bzeta^{i},\bchi^i )=\sum_{i=1}^{9} H_i.\label{EFNSED13}
\end{align}%
In obtaining \ref{EFNSED26}, we have made use of the following facts: (i) Since $\bchi^0 =\bO_{fh}\bu^0-\bP_h\bu^0=0$, we have
\begin{align*}
\sum_{i=1}^{n} \Delta t e^{2\alpha t_i} \partial_t\|\bchi^i\|^2_{L^2(\Omega_f)} & = e^{2\alpha t_n} \|\bchi^n\|^2_{L^2(\Omega_f)}
 - \sum_{i=1}^{n-1} (e^{2\alpha \Delta t}-1)  e^{2\alpha t_i} \|\bchi^i\|^2_{L^2(\Omega_f)}.
\end{align*}%
(ii) Since $\xi^0=O_{ph}\phi^0-L_h\phi^0=0$, a similar break for $\partial_t\|\xi^i\|^2_{L^2(\Omega_f)}$ as well. (iii) We have applied the fact that $e^{2\alpha \Delta t}-1\leq C(\alpha)\Delta t$. \\
\noindent We now begin to bound the trilinear terms $H_i$'s. Using (\ref{T9}), (\ref{nonliearity1}), and (\ref{aprioricondition1}), the first term can be bounded
\begin{align*}
\Delta t \sum_{i=1}^{n}  e^{2\alpha t_i}|H_1|\leq C_0\Delta t \sum_{i=1}^{n}  e^{2\alpha t_i}\|\bU_h^{i-1}\|_{\bX_1}\|\bchi^i\|^2_{\bX_1}\leq \frac{\mu C_3}{4}\Delta t \sum_{i=1}^{n}  e^{2\alpha t_i}\|\bchi^i\|^2_{\bX_1}.
\end{align*}
\noindent Young’s inequality and estimate (\ref{nonlinearity3}) ensure the next bound
\begin{align*}
\Delta t \sum_{i=1}^{n}  e^{2\alpha t_i}|H_2| 
\leq  \frac{\mu C_3}{40}\Delta t \sum_{i=1}^{n}  e^{2\alpha t_i}\|\bchi^{i}\|^2_{\bX_h} + C\Delta t \sum_{i=1}^{n}  e^{2\alpha t_i} \|\bu^i\|^2_{2,\Omega_f}\|\bchi^{i-1}\|^2_{L^2(\Omega_f)}.
\end{align*}
By means of inverse inequality (\ref{T4}), estimates (\ref{nonlinearity3}), (\ref{nonlinearity5}) and Young’s inequality, we obtain	
\begin{align*}
\Delta t \sum_{i=1}^{n}  e^{2\alpha t_i}|H_3+H_4|\leq  \frac{\mu C_3}{40}\Delta t \sum_{i=1}^{n}  e^{2\alpha t_i}\|\bchi^{i}\|^2_{\bX_1} + C \Delta t \sum_{i=1}^{n}  e^{2\alpha t_i}\bigg(\|\bu\|_{2,\Omega_f}\|\bchi^{i-1}\|^2_{L^2(\Omega_f)}+  \|\bzeta^i\|^2_{\bX_1}\|\bzeta^{i-1}\|^2_{\bX_1}\bigg).
\end{align*}%
To bound the next terms, we use
$\bu^i-\bu^{i-1}= \int_{t_{i-1}}^{t_i} \bu_t(t)\, dt$
and standard inequalities, to arrive at
\begin{align*}
\Delta t \sum_{i=1}^{n}  e^{2\alpha t_i}|H_5+H_6+H_7|
&\leq  \frac{\mu C_3}{40}\Delta t \sum_{i=1}^{n}  e^{2\alpha t_i}\|\bchi^{i}\|^2_{\bX_1} + C \Delta t^2 \sum_{i=1}^{n}  e^{2\alpha t_i}\bigg(\|\bu\|^2_{2,\Omega_f}\int_{t_{i-1}}^{t_i} \|\bu_t\|^2_{L^2(\Omega_f)} dt\\
&\quad\quad\quad\quad\quad+\|\bu^i\|^2_{2,\Omega_f}\int_{t_{i-1}}^{t_i} \|\bzeta_t\|_{\bX_1}^2 dt+ \|\bzeta^i\|^2_{\bX_1}\int_{t_{i-1}}^{t_i} \|\bu_t\|^2_{1,\Omega_f} dt\bigg).
\end{align*}
\noindent
Applying the standard inequalities along with (\ref{T8}) and (\ref{nonlinearity7}), we find from $H_8$ and $H_9$ that 
\begin{align*}
\Delta t \sum_{i=1}^{n}  & e^{2\alpha t_i}\Big| l_{h}(\bu^{i}, \bu^i; \bzeta^{i},\bchi^i )
 - l_{h}(\bU_h^{i-1}, \bu^i; \bzeta^{i},\bchi^i )  \Big|
\leq  \frac{\mu C_3}{40}\Delta t \sum_{i=1}^{n}  e^{2\alpha t_i}\|\bchi^{i}\|^2_{\bX_1}\nonumber\\
&\qquad +C\Delta t \sum_{i=1}^{n}  e^{2\alpha t_i}\bigg(\|\bzeta^{i-1}\|^2_{\bX_1}\|\bzeta^i\|^2_{\bX_1}     +\Delta t\|\bzeta^i\|^2_{\bX_1}\int_{t_{i-1}}^{t_i} \|\bu_t\|^2_{L^2(\Omega_f)} dt
+  \|\bu^i\|_{2,\Omega_f}  \|\bchi^{i-1}\|^2_{L^2(\Omega_f)}\bigg). 
\end{align*}
\noindent
We observe that $\|\bu^i\|^2_{2,\Omega_f}\leq C$ and $\|\bzeta\|_{\bX_1}$, $\|\bzeta_t\|_{\bX_1}\leq C$ from (\ref{Oseenprojection1}).
Also, terms appearing in $H_5,H_6, H_7$, etc. can be estimated, using assumption $(\textbf{A}\textbf{1})$, as follows: 
\begin{align*}
C {\Delta t}^2 \sum_{i=1}^{n}  e^{2\alpha t_i}\|\bu^i\|^2_{2,\Omega_f} \int_{t_{i-1}}^{t_i} \|\bu_t\|^2_{L^2(\Omega_f)} dt &\le C {\Delta t}^2 \sum_{i=1}^{n}   \int_{t_{i-1}}^{t_i} e^{2\alpha (t_i-s)} e^{2\alpha s}\|\bu_{t}\|^2_{L^2(\Omega_f)} dt \\
&\quad\leq C {\Delta t}^2 e^{2\alpha \Delta t} \int_{0}^{t_n}  e^{2\alpha s}\|\bu_{t}\|^2_{L^2(\Omega_f)} dt \leq C {\Delta t}^2 e^{2\alpha t_{n+1} }, \\
C{\Delta t}^2 \sum_{i=1}^{n}  e^{2\alpha t_i} \|\bu^i\|^2_{2,\Omega_f}\int_{t_{i-1}}^{t_i} \|\bzeta_t\|_{\bX_1}^2 dt &\leq  C {\Delta t}^2 \sum_{i=1}^{n}e^{2\alpha t_i} \int_{t_{i-1}}^{t_i} dt  \leq C\Delta t^3 \sum_{i=1}^{n}  e^{2\alpha t_i}\\
&\quad \leq C\Delta t^3 e^{2\alpha\Delta t}\frac{(e^{2\alpha n\Delta t}-1)}{e^{2\alpha\Delta t}-1}\leq C\Delta t^2 e^{2\alpha t_{n+1}}.
\end{align*}
Also, we have
\begin{align*}
C\Delta t \sum_{i=1}^{n}  e^{2\alpha t_i} \|\bzeta^{i-1}\|^2_{\bX_1}\|\bzeta^i\|^2_{\bX_1} 
\leq Ch^{4k_1+4k_2} \Delta t e^{2\alpha\Delta t}\dfrac{(e^{2\alpha n\Delta t}-1)}{(e^{2\alpha\Delta t}-1)}
 \leq Ch^{4k_1+4k_2} e^{2\alpha t_{n+1}}.
\end{align*}%
\noindent 
Observe that we have avoided the integration by parts in estimating the nonlinear terms above (see Remark \ref{bypartsresult} below). And that has been made possible by the application of modified Oseen operators, which have introduced time step difference, thus allowing the estimates in terms of $\Delta t$. 

\noindent Finally, collecting the nonlinear bounds, applying the discrete Gr\"{o}nwall's inequality, multiplying the resulting equation by $e^{-2\alpha t_n}$, and using triangle inequality, approximation properties (\ref{Oseenprojection1}) along with assumption $(\textbf{A}\textbf{1})$, the desired result is achieved. This concludes the proof. 
\end{proof}%
\noindent The next theorem is one of our key contributions, which states optimal $L^\infty(\bL^2)$ and $L^\infty(L^2)$-norm results for velocity and piezometric head errors, respectively, and it results as a combination of Lemmas \ref{Opproperty} and \ref{chixi}.
\begin{theorem} \label{velo1}
Let the hypothesis of Lemma \ref{chixi} be satisfied. Then there exists a constant $K_T>0$ such that 
 \begin{align}
\|\bu^i-\bU_h^i\|^2_{L^2(\Omega_f)} +  gS_0\|\phi^i-\Phi_h^i\|^2_{L^2(\Omega_p)}\leq K_T(h^{2k_1+2}+h^{2k_2+2}+\Delta t^2). 
\end{align}
\end{theorem}
\begin{Remark}\label{bypartsresult}
As discussed in the 'Introduction', the modified Oseen operators play a crucial role in establishing the above theorem, Theorem \ref{velo1}. A standard approach to derive $L^\infty(\bL^2)$ estimates for NSE is to perform integration by parts to handle the nonlinear term. However, in the case of NSD, if we calculate $c_{NS}(\bU_h^{i-1}, \bu^i;\bzeta^i,\bchi^i)$ using integration by parts, then from (\ref{trilinear}) we find
\begin{align*}
c_{NS}(\bU_h^{i-1}, \bu^i;\bzeta^i,\bchi^i)&=\sum_{E\in\mathcal{E}_{fh}}  (\bu^i\cdot \nabla \bchi^i,\bzeta^i)_E -\frac{1}{2}\sum_{E\in \mathcal{E}_{1h}}(\nabla\cdot \bu^i,\bzeta^i\cdot \bchi^i)  +\frac{1}{2}\sum_{e\in \Gamma_{fh}\cup\Gamma_f}([\bu^i]\cdot \bn_e,\{\bzeta^i \cdot \bchi^i\})\nonumber\\
& -  \sum_{E\in\mathcal{E}_{fh}}(|\bu^i\cdot\bn_{E}|\bzeta^{i,ext}, [\bchi^i])_{\partial E_-(\bU_h^{i-1}) \setminus \Gamma}+\frac{1}{2} \sum_{e\in \Gamma}(|\bu^i\cdot\bn_e|,\bzeta^i \cdot \bchi^i).
\end{align*}
The interface term can be estimated by means of (\ref{T1}), (\ref{T3}), (\ref{T8}), and (\ref{Oseenprojection1}) in the following three ways: 
\begin{align*}
| \sum_{e\in \Gamma}(|\bu^i\cdot\bn_e|,\bzeta^i \cdot \bchi^i)|&\leq \sum_{e\in \Gamma}\|\bu^i\|_{L^4(e)}\|\bzeta^i\|_{L^2(e)}\|\bchi^i\|_{L^4(e)}\\
&\leq C|\bu^i|_{1,\Omega_f}(h^{k_1+1/2}(|\bu|_{k_1+1,\Omega_f}
+|p|_{k_1,\Omega_f})+h^{k_2+1/2}|\phi|_{k_2+1,\Omega_p})\|\bchi^i\|_{\bX_1},\\
| \sum_{e\in \Gamma}(|\bu^i\cdot\bn_e|,\bzeta^i \cdot \bchi^i)|&\leq \sum_{e\in \Gamma}|e|^{1/4}\|\bu^i\|_{L^\infty(e)}\|\bzeta^i\|_{L^2(e)}\|\bchi^i\|_{L^4(e)}\\
&\leq C\|\bu^i\|_{L^\infty(\Omega_f)}(h^{k_1+3/4}(|\bu|_{k_1+1,\Omega_f}
+|p|_{k_1,\Omega_f})+h^{k_2+3/4}|\phi|_{k_2+1,\Omega_p})\|\bchi^i\|_{\bX_1},\\
| \sum_{e\in \Gamma}(|\bu^i\cdot\bn_e|,\bzeta^i \cdot \bchi^i)|&\leq \sum_{e\in \Gamma}|e|^{1/2}\|\bu^i\|_{L^\infty(e)}\|\bzeta^i\|_{L^4(e)}\|\bchi^i\|_{L^4(e)}\\
&\leq C\|\bu^i\|_{L^\infty(\Omega_f)}(h^{k_1+3/4}(|\bu|_{k_1+1,\Omega_f}
+|p|_{k_1,\Omega_f})+h^{k_2+3/4}|\phi|_{k_2+1,\Omega_p})\|\bchi^i\|_{\bX_1}.
\end{align*}
\noindent This shows that the best we can do is $O(h^{k_1+3/4})$, which is suboptimal.
\end{Remark}
\begin{Remark}\label{Oseenhelp}
As we have noticed in Lemma \ref{chixi}, the modified Oseen operators absorb all the interface terms in their definition. Consequently, they do not appear in the error equation. Therefore, the difficulties that arise regarding them, as discussed in 'Introduction', have been taken care of. Apart from that, we have noticed in the above lemma that integration by parts has been avoided by means of the modified Oseen operators. Thus, instead of suboptimal estimates in $h$, we have taken contributions in terms of $\Delta t$ from the last five trilinear terms in (\ref{EFNSED13}) and managed to obtain optimal results. 
\end{Remark}\vspace*{-2mm}
\subsection{Error Analysis for Pressure } \label{s7} 
Our final aim is to derive the optimal fully discrete pressure error bound in the $L^\infty(L^2)$-norm, which requires the estimates of $\partial_t\e_f^i$ and $\partial_t e_p^i$. These estimates, first in $L^2(\bL^2)$ and $L^2(L^2)$-norms, and then in $L^\infty(\bL^2)$ and $L^\infty(L^2)$-norms, are presented in two auxiliary lemmas below. The main technical challenge here is that, unlike earlier, the interface terms of the form $(\cdot,\cdot)_{\Gamma}$ do not get cancelled out. We have used (\ref{T12}) and summation by parts to handle the interface terms.
\begin{lemma}{\label{Intermidiate1}}
Let the hypothesis of Lemma \ref{chixi} be satisfied. Then the errors satisfy
\begin{align}
\|\e_f^n\|^2_{\bX_1} + \|e_p^n\|^2_{M_{1p}}+\Delta t e^{-2\alpha t_n} \sum_{i=1}^{n} e^{2\alpha t_i} \|\partial_t\e_f^n\|^2_{L^2(\Omega_f)} 
 + \Delta t e^{-2\alpha t_n} \sum_{i=1}^{n} e^{2\alpha t_i}&\|\partial_t e_p^n\|^2_{L^2(\Omega_p)} \nonumber\\
+ \Delta t e^{-2\alpha t_n} \sum_{i=1}^{n} e^{2\alpha t_i}\frac{1}{{G}}\|\partial_t\e_f^i\cdot\btau^j_{\gamma}\|^2_{L^2(\Gamma)} +\frac{1}{{G}}\|\e_f^n\cdot\btau_{\gamma}\|^2_{L^2(\Gamma)}& \leq  K_T(h^{2k_1}+h^{2k_2}+\Delta t^2).\label{errorestimate2}
\end{align}
\end{lemma}
\begin{proof}
First, we apply the definition of the modified Oseen operators (\ref{Oseenprojection}) in (\ref{EFNSED10}) and plug $\bv_h=\partial_t \bchi^i$, $\psi_h=\partial_t \xi^i$. As a consequence, we will observe that there will be no interface terms on the right hand side of the error equation. The bilinear terms on the right hand side of the resulting equation can be bounded following similar steps as in Lemma \ref{chixi}.

\noindent The nonlinear terms are rewritten as follows:
\begin{align}
c_{NS}&(\bU_h^{i-1},  \bU_h^{i-1};\bU_h^{i},\partial_t \bchi^i )
-c_{NS}(\bu^i,  \bu^i;\bu^i,\partial_t \bchi^i)
+c_{NS}(\bu^i,  \bu^i;\bzeta^i,\partial_t \bchi^i) 
 +c_{NS}(\bu^i,  \bzeta^i;\bu^i,\partial_t \bchi^i) \nonumber\\
&  = c_{NS}(\bU_h^{i-1},  \bchi^{i-1};\bzeta^{i},\partial_t \bchi^i )+ c_{NS}(\bU_h^{i-1},  \e_f^{i-1};\bu^i,\partial_t \bchi^i )+ c_{NS}(\bU_h^{i-1},  \bu^{i-1}-\bu^i;\bu^i,\partial_t \bchi^i)\nonumber\\
&\hspace*{3mm}+ c_{NS}(\bU_h^{i-1}, \bzeta^{i-1} ;\e_f^{i-1},\partial_t \bchi^i )- c_{NS}(\bU_h^{i-1},  \bu^{i-1};\e_f^{i},\partial_t \bchi^i) 
 +c_{NS}(\bu^i,  \bu^i;\bzeta^i,\partial_t \bchi^i) \nonumber\\
&\hspace*{3mm} +c_{NS}(\bu^i,  \bzeta^i;\bu^i,\partial_t \bchi^i)
 + c_{NS}(\bU_h^{i-1},  \bchi^{i-1};\bchi^{i},\partial_t \bchi^i ).\label{IEFNSED3} 
\end{align}%
The nonlinear terms can be bounded as in the previous Lemma \ref{chixi} except for the last term in (\ref{IEFNSED3}). Using (\ref{T6}) along with standard inequalities, the following term can be  estimated as
\begin{align*}
|c_{NS}(\bU_h^{i-1},  \bchi^{i-1};\bchi^{i},\partial_t \bchi^i )|\leq \frac{1}{28} \|\partial_t \bchi^i\|^2_{L^2(\Omega_f)}+ \frac{1}{\min\limits_{{E\in\mathcal{E}_{fh}}}h_E } \|\bchi^{i-1}\|^2_{\bX_h}\|\bchi^{i}\|^2_{\bX_h} . 
\end{align*}
\noindent 
Inserting all the bounds, multiplying by $\Delta t e^{2\alpha t_i}$ and summing from $i=1$ to $n \,(\leq N)$ we arrive at
\begin{align}
&\Delta t\sum_{i=1}^{n} e^{2\alpha t_i}\Bigg(\|\partial_t \bchi^i\|^2_{L^2(\Omega_f)} + gS_0\|\partial_t \xi^i\|^2_{L^2(\Omega_p)}+ a_{NS}(\bchi^i,\partial_t \bchi^i)+ ga_D(\xi^i,\partial_t \xi^i)+  \frac{1}{G}(\bchi^i\cdot\btau_{\gamma},\partial_t \bchi^i\cdot\btau_{\gamma})_{{\Gamma}} \Bigg)\nonumber\\
&\,\, +  g\Delta t\sum_{i=1}^{n} e^{2\alpha t_i}\Bigg((\xi^i,\partial_t \bchi^i\cdot{\bn_{\gamma}})_{{\Gamma}} 
-  (\bchi^i\cdot\bn_{\gamma},\partial_t \xi^i)_{{\Gamma}} \Bigg) 
\leq \frac{12\Delta t}{28}\sum_{i=1}^{n} e^{2\alpha t_i}  \|\partial_t \bchi^i\|^2_{L^2(\Omega_f)}\nonumber\\
 &\quad  + \frac{2gS_0\Delta t}{4}\sum_{i=1}^{n} e^{2\alpha t_i}  \|\partial_t \xi^i\|^2_{L^2(\Omega_f)} +C \frac{\Delta t}{\min\limits_{{E\in\mathcal{E}_{fh}}}h_E } \sum_{i=1}^{n} e^{2\alpha t_i} \|\bchi^{i-1}\|^2_{\bX_1}\|\bchi^{i}\|^2_{\bX_1}
  + C \Delta t \sum_{i=1}^{n} e^{2\alpha t_i}\bigg( \|\bchi^{i-1}\|^2_{\bX_1}\nonumber\\
 &\quad\quad + \|\e_f^{i-1}\|^2_{\bX_1} + \|\e_f^{i}\|^2_{\bX_1} + \|\bzeta^i\|^2_{\bX_h} \bigg)+C(h^{2k_1+2}+h^{2k_2+2} + 
 {\Delta t}^2 ).\label{IEFNSED2.O}
\end{align}
\noindent Note that the last few terms involving $h$ and $\Delta t$ in the above inequality are obtained from applying the regularity assumption ({\bf A1}) of the solutions $\bu$ and $\phi$. 

\noindent
For the terms on the left side of (\ref{IEFNSED2.O}), we shift the discrete time derivative to the other arguments in the bilinear forms by using summation by parts. Since $a_{NS}(.,.)$ is symmetric, we obtain
\begin{align*}
\Delta t \sum_{i=1}^{n}  e^{2\alpha t_i}  a_{NS}(\bchi^i,\partial_t \bchi^i) & = \frac{1}{2}e^{2\alpha t_n} a_{NS}(\bchi^n, \bchi^n)  - \frac{1}{2}  a_{NS}(\bchi^0,  \bchi^0)\nonumber\\
&\quad - \frac{1}{2}  \sum_{i=0}^{n-1} e^{2\alpha t_i} (e^{2\alpha \Delta t}-1)    a_{NS}(\bchi^i, \bchi^i)   + \frac{{\Delta t}^2}{2}  \sum_{i=1}^{n} e^{2\alpha t_i}  a_{NS}(\partial_t\bchi^i, \partial_t\bchi^i).
\end{align*}
Similar calculations can be applied to will go for $a_{D}(\cdot,\cdot)$. In the case of interface terms, we follow the same approach. For instance, the fifth term on the left side of (\ref{IEFNSED2.O}) can be handled as
\begin{align*}
&\Delta t \sum_{i=1}^{n}   \frac{1}{G} e^{2\alpha t_i}(\bchi^i\cdot\btau_{\gamma},\partial_t \bchi^i\cdot\btau_{\gamma})_{{\Gamma}}  = \frac{1}{2G} (e^{2\alpha t_n}(\bchi^n\cdot\btau_{\gamma}, \bchi^n\cdot\btau_{\gamma})_{{\Gamma}} -  (\bchi^0\cdot\btau_{\gamma}, \bchi^0\cdot\btau_{\gamma})_{{\Gamma}}) \nonumber\\
&\quad - \frac{1}{2G} \sum_{i=0}^{n-1} 
  e^{2\alpha t_i} (e^{2\alpha \Delta t}-1)   (\bchi^i\cdot\btau_{\gamma},\bchi^i\cdot\btau_{\gamma})_{{\Gamma}}   + \frac{{\Delta t}^2}{2G}  \sum_{i=1}^{n} e^{2\alpha t_i} (\partial_t\bchi^i\cdot\btau_{\gamma},\partial_t \bchi^i\cdot\btau_{\gamma})_{{\Gamma}} 
\end{align*}
To bound the sixth term on the left side of (\ref{IEFNSED2.O}), we employ (\ref{T12}), (\ref{erroresult1}) and Young's inequality as follows:
\begin{align*}
\Delta t \sum_{i=1}^{n}  e^{2\alpha t_i} \big|   g(\xi^i,\partial_t \bchi^i\cdot{\bn_{\gamma}})_{{\Gamma}} \big| &\leq \frac{\Delta t}{28}  \sum_{i=1}^{n}  e^{2\alpha t_i} \|\partial_t \bchi^i\|^2_{L^2(\Omega_f)}  + C(h^{2k_1+2}+h^{2k_2+2}+{\Delta t}^2). 
\end{align*}
However, to bound the seventh term, summation by parts has been used, since (\ref{T12}) is no longer applicable.
\begin{align*}
&\Delta t \sum_{i=1}^{n}  e^{2\alpha t_i} \big|  g(\bchi^i\cdot\bn_{\gamma},\partial_t \xi^i)_{{\Gamma}} \big| \leq g e^{2\alpha t_n}\big|(\bchi^n\cdot\bn_{\gamma}, \xi^n)_{{\Gamma}}\big| + ge^{2\alpha t_0}\big|(\bchi^0\cdot\bn_{\gamma}, \xi^0)_{{\Gamma}} \big|\nonumber\\
&\quad \quad + g\sum_{i=1}^{n} e^{2\alpha t_{i-1}}(e^{2\alpha \Delta t}-1) \big|(\bchi^{i-1}\cdot\bn_{\gamma}, \xi^{i-1})_{{\Gamma}}\big| + g\Delta t \sum_{i=1}^{n}e^{2\alpha t_i} \big|(\partial_t\bchi^i\cdot\bn_{\gamma}, \xi^{i-1})_{{\Gamma}} \big|\nonumber\\
&\qquad\quad \leq    \frac{\Delta t}{28}\sum_{i=1}^{n}  e^{2\alpha t_i} \|\partial_t \bchi^i\|^2_{L^2(\Omega_f)} + \frac{gC_4}{2} e^{2\alpha t_n}\|\xi^n  \|^2_{M_{1p}} + C e^{2\alpha t_n} \frac{1}{2} \|\bchi^n  \|^2_{L^2(\Omega_f)}\nonumber\\
 &\qquad\quad\quad\,\,\,\, + C \Delta t \sum_{i=1}^{n} ( e^{2\alpha t_{i-1}}+e^{2\alpha t_i}) \|\xi^{i-1}  \|^2_{M_{1p}} + C\Delta t \sum_{i=1}^{n}  e^{2\alpha t_i} \|\bchi^{i-1}  \|^2_{\bX_1}.
 \end{align*}
\noindent We insert all the above bounds in (\ref{IEFNSED2.O}) and utilize Lemma \ref{coercivity}. Then a use of the discrete Gr\"{o}nwall inequality, (\ref{erroresult1}), Theorem \ref{velo1} and $\textbf{(A1)}$ yields
\begin{align}
\Delta t  \sum_{i=1}^{n}  e^{2\alpha t_i}&\bigg(\|\partial_t \bchi^i\|^2_{L^2(\Omega_f)} + gS_0 \|\partial_t \xi^i\|^2_{L^2(\Omega_p)}+   \mu C_3 {\Delta t} \|\partial_t \bchi^i\|^2_{L^2(\Omega_f)} 
+  gC_4 {\Delta} t  \|\partial_t \xi^i\|^2_{L^2(\Omega_f)} \nonumber\\
&\quad  +  \frac{\Delta t}{{G}}\|\partial_t\bchi^{n}   \cdot\btau_{\gamma}\|^2_{L^2(\Gamma)}\bigg)+ e^{2\alpha t_n}\bigg( \mu C_3 \|\bchi^{n}  \|^2_{\bX_1}  +  gC_4  \|\xi^n  \|^2_{M_{1p}} 
  +    \frac{1}{{G}}\|\bchi^{n}  \cdot\btau^j_{\gamma}\|^2_{L^2(\Gamma)}\bigg) \nonumber\\
&\quad\quad\quad\leq C e^D(h^{2k_1}+h^{2k_2}+{\Delta t}^2),\label{IEFNSED13}
\end{align}
\noindent where $D=\frac{h^{2k_1}+{\Delta t}^2}{\min\limits_{{E\in\mathcal{E}_{fh}}}h_E }$. Here, we necessitate that the mesh fulfils $D\leq L$ for some $L>0$, which does not depend on $h$ and $\Delta t$. Finally, the use of triangular inequality and Lemma \ref{Oseenprojection1}, assumption (\textbf{A1}), and multiplication by $e^{-2\alpha t_n}$ provides our desired result.
\end{proof}
\begin{Remark}
We would like to mention that the restriction $D\leq L$ is milder than the quasi-uniformity assumption in space (please see \cite[(109)]{CGPR17}).
\end{Remark}
\begin{Remark}\label{pressure}
In Lemma \ref{chixi}, we have observed that the interface terms in the right hand side of the error equation are absorbed by modified Oseen operators, and one of the interface terms on the left hand side of that equation gives a norm, and the other two get cancelled out. However, in this lemma, the scenario is quite different, as the interface terms on the left hand side of the error equation do not cancel out. To bound the first one that is $(\xi^i,\partial_t \bchi^i\cdot{\bn_{\gamma}})_{{\Gamma}}$ we employ (\ref{T12}). But for the other one, i.e. $(\bchi^i\cdot\bn_{\gamma},\partial_t \xi^i)_{{\Gamma}}$, we can not follow a similar technique, as we need to bound $\|\partial_t \xi^i\|$ in $L^2$-norm, not in energy norm. Therefore, to overcome the difficulty, we have used summation by parts and achieved the desired result. 
\end{Remark}
\noindent For the next auxiliary lemma, we have used standard projections instead of modified Oseen operators for two reasons. First, we require only $O(h^k)$ estimates and not $O(h^{k+1})$; second, the use of modified Oseen operators would bring in extra nonlinear terms, adding to additional technicalities which can now be avoided.

\noindent Let us begin by defining the discretization velocity and piezometric head errors $\bLambda^i = \bP_h\bu^i-\bU_h^i$ and $\Upsilon^i = L_h\phi^i-\Phi_h^i$, respectively. We then rewrite (\ref{EFNSED9}) as follows:
\begin{align}
&(\partial_t \bLambda^i,\bv_h)+gS_0(\partial_t \Upsilon^i, \psi_h)+a_{NS}(\bLambda^i,\bv_h)+ ga_D(\Upsilon^i, \psi_h) + \frac{1}{G}(\bLambda^i\cdot\btau^j_{\gamma},\bv_h\cdot\btau^j_{\gamma})_{{\Gamma}}\nonumber\\
& \quad +c_{NS}(\bu^i,  \bu^i;\bu^i, \bv_h)
-c_{NS}(\bU_h^{i-1},  \bU_h^{i-1};\bU_h^{i},\bv_h )
+  g(\Upsilon^i,\bv_h\cdot{\bn_{\gamma}})_{{\Gamma}}  - g(\bLambda^i\cdot\bn_{\gamma}, \psi_h)_{{\Gamma}}\nonumber\\
 &\quad +  b_{NS}(\bv_h,p^i)=(\partial_t \bu^i - \bu_t^i,\bv_h) +  gS_0(\partial_t \phi^i- \,\,\phi_t^i, \psi_h) -(\partial_t(\bu^i-\bP_h\bu^i),\bv_h) -(\partial_t(\phi^i-L_h\phi^i),\psi_h) \nonumber\\
&\quad\quad
 -  a_{NS}(\bu^i-\bP_h\bu^i,\bv_h) - ga_D(\phi^i-L_h\phi^i, \psi_h) 
 -  g( \phi^i-L_h\phi^i,\bv_h\cdot{\bn_{\gamma}})_{{\Gamma}} 
+  g(\bu^i-\bP_h\bu^i\cdot\bn_{\gamma},\psi_h)_{{\Gamma}} \nonumber\\
&\quad\quad- \frac{1}{G}((\bu^i-\bP_h\bu^i)\cdot\btau^j_{\gamma},\psi_h\cdot\btau^j_{\gamma})_{{\Gamma}},  \quad \quad \forall \bv_h\in\bV_h,\psi_h\in M_{ph}. \label{AEFNSED1}
\end{align}
\noindent The following lemma establishes $L^{\infty}(\bL^2)$ and $L^{\infty}(L^2)$-norm estimates for $\partial_t \e_f^i$ and $\partial_t e_p^i$, respectively.
\begin{lemma}{\label{Intermidiate2}}
The errors $\partial_t \e_f^i$ and $\partial_t e_p^i$ satisfy
\begin{align}
\|\partial_t \e_f^n\|^2_{L^2(\Omega_f)} + \|\partial_t e_p^n\|^2_{L^2(\Omega_p)} + e^{-2\alpha t_n} \Delta t  &\sum_{i=1}^{n}  e^{2\alpha t_i}\|\partial_t \e_f^i\|^2_{\bX_1} + e^{-2\alpha t_n} \Delta t  \sum_{i=1}^{n}  e^{2\alpha t_i}\|\partial_t e_p^i\|^2_{M_{p}} \nonumber\\
&\leq K_T(h^{2k_1}+h^{2k_2}+\Delta t ) .\label{AEFNSED2}
\end{align}
\end{lemma}
\begin{proof}
Subtract (\ref{AEFNSED1}) at time $t=t_{i-1}$ from the one at time $t=t_{i}$, and substitute $\bv_h =\partial_t \bLambda^i $ and $\psi_h= \partial_t\Upsilon^i$. Then apply Lemma \ref{coercivity} on the left hand side of the resulting equation.

\noindent Now, let us rewrite the trilinear terms as follows
\begin{align}
&c_{NS}(\bU_h^{i-1},  \bU_h^{i-1};\bU_h^{i},\partial_t\bLambda^i ) 
-  c_{NS}(\bU_h^{i-2},  \bU_h^{i-2};\bU_h^{i-1},\partial_t\bLambda^i ) 
-c_{NS}(\bu^i,  \bu^i;\bu^i, \partial_t \bLambda^i) \nonumber\\
&\quad +  c_{NS}(\bu^{i-1},  \bu^{i-1};\bu^{i-1}, \partial_t \bLambda^i)\nonumber\\
&=- \Delta t c_{NS}(\bU_h^{i-1},  \partial_t\bu^{i-1};\e_f^{i},\partial_t\bLambda^i )- \Delta t c_{NS}(\bU_h^{i-1},  \partial_t (\bu^{i}- \bu^{i-1});\bu^{i},\partial_t\bLambda^i )\nonumber\\
&\quad- \Delta t c_{NS}(\bU_h^{i-1},  \partial_t\e_f^{i-1};\bU_h^{i},\partial_t\bLambda^i )- \Delta t c_{NS}(\bU_h^{i-1},   \bu^{i-1};\partial_t\e_f^{i},\partial_t\bLambda^i )\nonumber\\
&\quad- \Delta t c_{NS}(\bU_h^{i-1},  \bu^{i-1}- \bu^{i-2};\partial_t\bU_h^{i},\partial_t\bLambda^i )- \Delta t c_{NS}(\bU_h^{i-1},  \e_f^{i-2};\partial_t\bU_h^{i},\partial_t\bLambda^i )\nonumber\\
&\quad + l(\bU_h^{i-1},  \bU_h^{i-2};\bU_h^{i-1},\partial_t\bLambda^i )
- l(\bU_h^{i-2},  \bU_h^{i-2};\bU_h^{i-1},\partial_t\bLambda^i )=\sum_{i=1}^{8}M_i.
\label{AEFNSED4}
\end{align}
Here, the first two nonlinear terms can be bounded following similar calculations shown in previous lemmas. 
Since $\e_f^{i}=\bu^{i-1}-\bP_h\bu^{i-1}+\bLambda^{i-1}$, we can rewrite the third term  
\begin{align*}
M_3 &= \Delta t c_{NS}(\bU_h^{i-1}, \partial_t (\bu^{i-1}-\bP_h\bu^{i-1});\bU_h^{i},\partial_t\bLambda^i )+ \Delta t c_{NS}(\bU_h^{i-1}, \partial_t \bLambda^{i-1};\bU_h^{i},\partial_t\bLambda^i )=M_{3,1}+M_{3,2}.
\end{align*}
The technique of bounding $M_{3,1}$ is standard, whereas $M_{3,2}$ can be bounded using (\ref{T4}), (\ref{T8}), and (\ref{nonlinearity8}) as
\begin{align*}
|M_{3,2}|
&\leq C\Delta t \|\partial_t\bLambda^{i-1}\|_{\bX_1}\|\bU_h^{i}\|_{\bX_1}\|\partial_t\bLambda^i\|_{L^4(\Omega_f)}\nonumber\\
&\leq \frac{ \mu C_3\Delta t}{32}e^{-2\alpha\Delta t}\|\partial_t\bLambda^{i-1}\|^2_{\bX_1} + \frac{ \mu C_3\Delta t}{32}\|\partial_t\bLambda^{i}\|^2_{\bX_1} + C\Delta t e^{4\alpha\Delta t}\|\bU_h^{i}\|^4_{\bX_1}\|\partial_t\bLambda^{i}\|^2_{L^2(\Omega_f)}.
\end{align*}
\noindent For the next three terms, we can rewrite the terms and bound them using similar steps as above. The last two upwinding terms can be estimated using (\ref{nonlinearity8}), (\ref{nonlinearity7}), along with the standard inequalities.

\noindent Now, incorporate the above bounds, multiply the resulting inequality by $e^{2\alpha t_i}$ and sum over $2\leq i\leq n\leq N$. Also, we use the fact that $\|\partial_t \bLambda^{i}\|_{L^2(\Omega_f)} \leq \|\partial_t\e_f^i\|_{L^2(\Omega_f)} + \|\partial_t(\bu^i-\bP_h\bu^i)\|_{L^2(\Omega_f)}$ and $\|\bLambda^{i}\|_{\bX_h}\leq \|\e_f^i\|_{\bX_h}+ \|\bu^i-\bP_h\bu^i\|_{\bX_1}$. Similar inequalities hold for $\|\partial_t\Upsilon\|^2_{L^2(\Omega_p)}$ and $\|\Upsilon\|^2_{M_{1p}}$. Along with these bounds, using (\ref{Ph}), we arrive at
\begin{align}
&\Delta t\sum_{i=2}^n e^{2\alpha t_i}\Bigg((\partial_t(\partial_t \bLambda^i),\partial_t \bLambda^i )+  gS_0(\partial_t(\partial_t \Upsilon^i), \partial_t\Upsilon^i)+    \mu C_3\|\partial_t \bLambda^i \|^2_{\bX_1} 
+  gC_4  \|\partial_t\Upsilon^i \|^2_{M_{1p}} +   \frac{1}{{G}}\|\partial_t \bLambda^i  \cdot\btau_{\gamma}\|^2_{L^2(\Gamma)} \Bigg )\nonumber\\
&\quad
\leq  \Delta t\sum_{i=2}^n e^{2\alpha t_i}\Bigg( (\partial_t(\partial_t \bu^i   - \bu_t^i ),\partial_t \bLambda^i) +   gS_0  t(\partial_t(\partial_t \phi^i  - \phi_t^i ),\partial_t \Upsilon^i)-   (\partial_t (\bu^i-\bP_h\bu^i),\partial_t \bLambda^i  )\nonumber\\
&\quad\quad -   gS_0  (\partial_t (\phi^i-L_h\phi^i), \partial_t \Upsilon^i)-   a_{NS}(\partial_t (\bu^i-\bP_h\bu^i),\partial_t \bLambda^i  )-   gS_0  a_D(\partial_t (\phi^i-L_h\phi^i), \partial_t \Upsilon^i)-    b_{NS}(\partial_t \bLambda^i, \partial_t p^i)
\nonumber\\ 
&\quad\quad\quad+  g  (\partial_t (\bu^i-\bP_h\bu^i)\cdot\bn_{\gamma},\partial_t \Upsilon^i)_{\Gamma}   - g  (\partial_t (\phi^i-L_h\phi^i), \partial_t \bLambda^i\cdot\bn_{\gamma})_{\Gamma}-   \frac{1}{{G}} (\partial_t (\bu^i-\bP_h\bu^i)\cdot\btau_{\gamma},\partial_t \bLambda^i \cdot\btau_{\gamma})_{\Gamma}\Bigg)\nonumber\\ 
&\quad\quad\quad + \frac{14\Delta t}{32}\sum_{i=2}^n e^{2\alpha t_i}\|\partial_t\bLambda^{i}\|^2_{\bX_1} + \frac{\Delta t}{4}\sum_{i=2}^n e^{2\alpha t_i}\|\partial_t\Upsilon^{i}\|^2_{\bX_1} + C e^{2\alpha t_n}(h^{2k_1}+h^{2k_2}+\Delta t).\label{AEFNSED3.O}
\end{align}
Here, the terms involving derivatives of $\bu$ and $\phi$ have been bounded in various norms using regularity assumptions. Next, we will bound the terms on the right-hand side of (\ref{AEFNSED3.O}).
The use of standard inequalities, Lemmas \ref{TL1}, \ref{TL2}, and the approximation property (\ref{projection2}) implies
\begin{align*}
\Delta t \sum_{i=2}^{n} e^{2\alpha t_i}\Big| a_{NS}(\partial_t (\bu^i-\bP_h\bu^i),\partial_t \bLambda^i  ) \Big| & \leq  \frac{ \mu C_3\Delta t}{32} \sum_{i=2}^{n} e^{2\alpha t_i} \|\partial_t\bLambda^{i}\|^2_{\bX_1}
 +  C h^{2k_1} e^{2\alpha \Delta t} \int_{t_{1}}^{t_{n}} e^{2\alpha t} |\bu_t|^2_{k_1+1,\Omega_f} dt.
\end{align*}
Using similar techniques, the third, fourth, and sixth terms can also be bounded. Now, for the interface term, using (\ref{T1}), (\ref{T9}), and (\ref{projection2}) along with the standard inequalities, we can derive 
\begin{align*}
&\Delta t \sum_{i=2}^{n} e^{2\alpha t_i}\Big| (\partial_t (\bu^i-\bP_h\bu^i)\cdot\bn_{\gamma},\partial_t \Upsilon^i)_{\Gamma}\Big|
\leq \frac{C_4\Delta t}{8}\sum_{i=2}^{n} e^{2\alpha t_i}\|\partial_t \Upsilon^i\|^2_{M_{1p}}   + Ch^{2k_1+1}e^{2\alpha \Delta t }\int_{t_{1}}^{t_{n}} e^{2\alpha t} |\bu_t|^2_{k_1+1} dt.
\end{align*}
Similarly, the ninth and tenth terms can be bounded using (\ref{projection3}). An application of (\ref{projection1}), the definition of $\bV_h$, the Cauchy-Schwarz inequality and Young's inequality yields 
 \begin{align*}
 \Delta t \sum_{i=2}^{n} e^{2\alpha t_i} \Big| b_{NS}(\partial_t \bLambda^i, \partial_t p^i) \Big| &= \Delta t \sum_{i=2}^{n} e^{2\alpha t_i} b_{NS}(\partial_t \bLambda^i, \partial_t( p^i-r_hp^i)) \nonumber\\
&\quad\leq \frac{ \mu C_3\Delta t}{32}\sum_{i=2}^{n} e^{2\alpha t_i}\|\partial_t \bLambda^i\|^2_{\bX_h} + Ch^{2k_1}e^{2\alpha \Delta t }\int_{t_{1}}^{t_{n}}e^{2\alpha t} |p_t|^2_{k_1} dt.
 \end{align*}
Note that we can write the first term of (\ref{AEFNSED3.O}) as follows:
\begin{align*}
\Delta t (\partial_t(\bu_t^i  - \partial_t \bu^i ),\partial_t \bLambda^i) = \Delta t ((\bu_t^i  - \partial_t \bu^i ),\partial_t \bLambda^i) -   (\bu_t^{i-1}  - \partial_t \bu^{i-1} ,\partial_t \bLambda^i )).  
\end{align*}
Therefore, using summation by parts, we obtain
\begin{align*}
&\Delta t \sum_{i=2}^{n} e^{2\alpha t_i} (\partial_t(\bu_t^i  - \partial_t \bu^i ),\partial_t \bLambda^i)  =  e^{2\alpha t_n} ((\bu_t^n  - \partial_t \bu^n ),\partial_t \bLambda^n)
 -  e^{2\alpha t_1} ((\bu_t^1  - \partial_t \bu^1 ),\partial_t \bLambda^1) \nonumber\\
& \quad - (e^{2\alpha \Delta t}-1) \sum_{i=2}^{n}  e^{2\alpha t_{i-1}} ((\bu_t^{i-1}  - \partial_t \bu^{i-1} ),\partial_t \bLambda^{i-1})- \Delta t \sum_{i=2}^{n} e^{2\alpha t_i} ((\bu_t^{i-1}  - \partial_t \bu^{i-1} ),\partial_t (\partial_t\bLambda^{i-1}))
\end{align*}
The above terms can be bounded by applying (\ref{Taylor}) and the standard inequalities.
\noindent Similarly, the second term can also be estimated. Substituting the above bounds in (\ref{AEFNSED3.O}) and applying (\ref{projection1}), (\ref{projection2}), Lemma \ref{Intermidiate1} and assumption $(\textbf{A1})$, we arrive at   
 \begin{align}
 e^{2\alpha t_n}\|\partial_t \bLambda^n\|^2_{L^2(\Omega_f)}  &+ e^{2\alpha t_n}\|\partial_t\Upsilon^n\|^2_{L^2(\Omega_p)} + \Delta t \sum_{i=2}^{n} e^{2\alpha t_i}\Bigg(\|\partial_t \bLambda^i\|^2_{\bX_h}
  + \|\partial_t\Upsilon^i\|^2_{M_{ph}}
   +  \frac{1}{{G}}\|\partial_t \bLambda^i\cdot\btau_{\gamma}\|^2_{L^2(\Gamma)}\Bigg)\nonumber\\
&  \leq C e^{2\alpha t_1}\|\partial_t \bLambda^1\|^2_{L^2(\Omega_f)}+ e^{2\alpha t_1}\|\partial_t\Upsilon^1\|^2_{L^2(\Omega_p)} 
  + C e^{2\alpha t_n}(h^{2k_1}+h^{2k_2}+\Delta t).\label{AEFNSED19}
 \end{align}
\noindent The above bound is not valid for $i=1$. To bound $\partial_t \bLambda^1$ and $\partial_t\Upsilon^1$, use the definition of modified Oseen operators (\ref{Oseenprojection}) in (\ref{EFNSED10}) with $i=1$, $\bv_h=\partial_t\bchi^1$, $\psi_h=\partial_t \xi^1$ and rewrite the nonlinear terms as
\begin{align}
&\|\partial_t\bchi^1\|^2_{L^2(\Omega_f)}+ gS_0\|\partial_t \xi^1\|^2_{L^2(\Omega_p)}+a_{NS}(\bchi^1,\partial_t\bchi^1)+ ga_D(\xi^1, \partial_t \xi^1) + \frac{1}{G}((\bchi^1\cdot\btau_{\gamma}),\partial_t\bchi^1\cdot\btau_{\gamma})_{{\Gamma}}\nonumber\\
& +  g(\xi^1,\partial_t\bchi^1\cdot{\bn_{\gamma}})_{{\Gamma}}- g(\bchi^1\cdot\bn_{\gamma}, \partial_t \xi^1)_{{\Gamma}}=  (\partial_t \bu^1 - \bu_t^1,\partial_t\bchi^1) +  gS_0(\partial_t \phi^1- \,\,\phi_t^1, \partial_t \xi^1) \nonumber\\
 & \quad -(\partial_t\bzeta^1,\partial_t\bchi^1) -  g(\partial_t\eta^1,\partial_t \xi^1)  -c_{NS}(\bU_h^{0}, \bu^0 ;\e_f^{1},\partial_t\bchi^1 )
 -c_{NS}(\bU_h^{0}, \e_f^{0} ;\bu^1,\partial_t\bchi^1 )
  +c_{NS}(\bU_h^{0},  \e_f^{0};\e_f^{1},\partial_t\bchi^1 )
  \nonumber\\
  &\quad-c_{NS}(\bU_h^{0}, \bu^1-\bu^0 ;\bu^1,\partial_t\bchi^1 )
+c_{NS}(\bu^1, \bu^1 ;\bzeta^1,\partial_t\bchi^1)
+c_{NS}(\bu^1, \bzeta^1 ;\bu^1,\partial_t\bchi^1 ). \label{AEFNSED21}
\end{align}
We recall the fact that $\bchi^0=0$ and $\xi^0=0$.
Now to handle the third and fourth terms on the left hand side of (\ref{AEFNSED21}), we use Lemma \ref{coercivity} and obtain
\begin{align*}
&\mu C_3\Delta t\|\partial_t\bchi^1\|^2_{\bX_1}\leq a_{NS}(\bchi^1,\partial_t\bchi^1) = \Delta t  a_{NS}(\partial_t\bchi^1,\partial_t\bchi^1),\,\,C_4\Delta t\|\partial_t\xi^1\|^2_{M_{1p}}\leq a_{D}(\xi^1, \partial_t \xi^1) = \Delta t  a_{D}(\partial_t \xi^1,\partial_t \xi^1).
\end{align*}
For the next terms on the left hand side of {(\ref{AEFNSED21})}, we arrive at 
\begin{align*}
&\quad\frac{1}{G}(\bchi^1\cdot\btau_{\Gamma},\partial_t\bchi^1\cdot\btau_{\Gamma})_{{\Gamma}} = \frac{1}{G}((\bchi^1-\bchi^0)\cdot\btau_{\Gamma},\partial_t\bchi^1\cdot\btau_{\Gamma})_{{\Gamma}}= \frac{1}{G}(\partial_t\bchi^1\cdot\btau_{\Gamma},\partial_t\bchi^1\cdot\btau_{\Gamma})_{{\Gamma}},\\
&- g(\xi^1,\partial_t\bchi^1\cdot{\bn_{\Gamma}})_{{\Gamma}} + g(\bchi^1\cdot\bn_{\Gamma},\partial_t\xi^1)_{{\Gamma}}
= - g\Delta t(\partial_t\xi^1,\partial_t\bchi\cdot{\bn_{\Gamma}})_{{\Gamma}} + g\Delta t(\partial_t\bchi^1\cdot\bn_{\Gamma},\partial_t\xi^1)_{{\Gamma}}=0.
\end{align*}
The rest of the terms on the right-hand side of (\ref{AEFNSED21}) can be estimated as Lemma \ref{Intermidiate1}. 
Thus, a combination of the above bounds, a multiplication of the resulting inequality by $e^{2\alpha t_1}$, and an application of (\ref{Oseenprojection2}), Lemma \ref{Intermidiate1}, and assumption $\textbf{(A1)}$ imply
\begin{align}
e^{2\alpha t_1}\bigg( \|\partial_t\bchi^1\|_{L^2(\Omega_f)}  +gS_0\|\partial_t \xi^1\|^2_{L^2(\Omega_p)} + \mu C_3\Delta t\|\partial_t\bchi^1\|^2_{\bX_1} &+ C_4\Delta t\|\partial_t\xi^1\|^2_{M_{1p}}
 +   \frac{\Delta t}{{G}}\|\partial_t\bchi^1\cdot\btau_{\gamma}\|^2_{{\Gamma}}\bigg)\nonumber\\
& \leq C(h^{2k_1}+ h^{2k_1}+\Delta t ).\label{AEFNSED26}
\end{align}
We now use the above estimates and the triangle inequality to obtain \vspace*{-2mm}
\begin{align*}
\|\partial_t\bLambda^1\|_{L^2(\Omega_f)} &\leq \|\partial_t(\bu^1-\bP_h\bu^1)\|_{L^2(\Omega_f)}+ \|\partial_t\bzeta^1\|_{L^2(\Omega_f)} + \|\partial_t\bchi^1\|_{L^2(\Omega_f)},\\
\|\partial_t\bLambda^1\|_{\bX_1} &\leq \|\partial_t(\bu^1-\bP_h\bu^1)\|_{\bX_1}+ \|\partial_t\bzeta^1\|_{\bX_1} + \|\partial_t\bchi^1\|_{\bX_1}.
\end{align*}
Similar estimates can be derived for $\|\partial_t\Upsilon^1\|_{L^2(\Omega_f)}$ and $\|\partial_t\Upsilon^1\|_{M_{1p}}$.

\noindent Next, using approximation properties (\ref{projection3}), (\ref{projection2}), and (\ref{Oseenprojection2}), we find  
\begin{align}
e^{2\alpha t_1}\bigg( \|\partial_t\bLambda^1\|_{L^2(\Omega_f)}  + gS_0\|\partial_t \Upsilon^1\|^2_{L^2(\Omega_p)}
 + \mu C_3\Delta t\|\partial_t\bLambda^1\|^2_{\bX_1} &+ C_4\Delta t\|\partial_t\Upsilon^1\|^2_{M_{1p}}
 +   \frac{\Delta t}{{G}}\|\partial_t\bLambda^1\cdot\btau_{\gamma}\|^2_{{\Gamma}}\bigg)\nonumber\\
& \leq C(h^{2k_1}+ h^{2k_1}+\Delta t ).\label{AEFNSED27}
\end{align}
Finally, incorporating (\ref{AEFNSED27}) in (\ref{AEFNSED19})
and using the relations $\partial_t\e_f^i=\partial_t(\bu^i-\bP_h\bu^i)+\partial_t\bLambda^i$, $\partial_t e_p^i=\partial_t(\phi^i-L_h\phi^i)+\partial_t\Upsilon^i$, (\ref{projection3}), (\ref{projection2}), and applying assumption $\textbf{(A1)}$, the remaining part of the proof is concluded.
\end{proof}
\noindent The next theorem presents the main result of the fully discrete $L^\infty(L^2)$-norm estimate for the pressure error.
\begin{theorem}
Assuming the hypothesis of Lemma \ref{chixi}, the following bound holds:
\begin{align}\label{pressurefinal}
\|p^i-P_h^i\|_{L^2(\Omega_f)} \leq K_T(h^{k_1}+h^{k_2}+{\Delta t}^{1/2} ).
\end{align}
\end{theorem}
\begin{proof}
Due to consistency of the numerical scheme, we observe that the weak solution $(\bu,p,\phi)$ satisfies (\ref{FNSED1}) with $\partial_t\bU_h^i$ replaced by $\bu_t^i$ at $t=t_i$.   
Therefore, subtracting (\ref{FNSED1}) from the resulting equation and choosing test functions $(\bv_h,0)$, we obtain the error equation as follows:
\begin{align}
&-b_{NS}(\bv_h,r_hp^i-P_h^i) =  (\bu_t^i - \partial_t \bu^i,\bv_h) + (\partial_t \e_f^i,\bv_h) + a_{NS}( \e_f^i,\bv_h)  + \frac{1}{G}(\e_f^i\cdot\btau_{\gamma},\bv_h\cdot\btau_{\gamma})_{{\Gamma}}\nonumber\\
 &\quad +  g(e_p^i,\bv_h\cdot\bn_{\gamma})_{{\Gamma}} +c_{NS}(\bu^i,  \bu^i;\bu^i, \bv_h)
-c_{NS}(\bU_h^{i-1},  \bU_h^{i-1};\bU_h^{i},\bv_h ) + b_{NS}(\bv_h,p^i-r_hp^i).\label{P1}
\end{align}
Here, the fifth term can be bounded using the standard inequalities with (\ref{T9}) and Lemma \ref{Intermidiate1} as
\begin{align}
\Big| \frac{1}{G}(\e_f^i\cdot\btau_{\gamma},\bv_h\cdot\btau_{\gamma})_{{\Gamma}}\Big| &\leq \frac{1}{{G}}\|\e_f^i\cdot\btau_{\gamma}\|_{L^2(\Gamma)}\|\bv_h\cdot\btau_{\gamma}\|_{L^2(\Gamma)}\leq C \frac{1}{{G}}\|\e_f^i\cdot\btau_{\gamma}\|_{L^2(\Gamma)}\|\bv_h\|_{\bX_1}\nonumber\\
&\quad\quad\leq  K_T(h^{k_1}+h^{k_2}+\Delta t)\|\bv_h\|_{\bX_1}.\label{P2}
\end{align}  
For the sixth term, using the standard inequalities, (\ref{T1}), (\ref{T9}), (\ref{Oseenprojection1}), and Lemma \ref{Intermidiate1}, we {acquire} 
\begin{align}
\Big| (e_p^i,\bv_h\cdot\bn_{\gamma})_{{\Gamma}} \Big| &\leq C \|e_p^i\|_{L^2(\Gamma)}\|\bv_h\|_{L^2(\Gamma)}
\leq  K_T(h^{k_1}+h^{k_2}+\Delta t)\|\bv_h\|_{\bX_1}.\label{P3}
\end{align}
Now, rewriting the trilinear terms, applying the continuity of $\bu^i$ and using (\ref{nonlinearity5}), we find
\begin{align}
&c_{NS}(\bu^i,  \bu^i;\bu^i, \bv_h)
- c_{NS}(\bU_h^{i-1},  \bU_h^{i-1};\bU_h^{i},\bv_h ) = c_{NS}(\bU_h^{i-1},  \bu^i;\bu^i, \bv_h) -c_{NS}(\bU_h^{i-1},  \bU_h^{i-1};\bU_h^{i},\bv_h ) \nonumber\\
&\quad= c_{NS}(\bU_h^{i-1},  \bu^i;\e_f^i, \bv_h)+ c_{NS}(\bU_h^{i-1},  \bu^i-\bu^{i-1};\bU_h^i, \bv_h)+ c_{NS}(\bU_h^{i-1},  \e_f^{i-1};\bU_h^i, \bv_h)\nonumber\\
&\quad\leq C (\|\bu^i\|_{1,\Omega_f}\|\e_f^i\|_{\bX_1} + \Delta t\|\bu_t\|_{L^{\infty}(0,T,H^1(\Omega_f))}\|\bU_h^{i}\|_{\bX_1}+ \|\e_f^{i-1}\|_{\bX_1}\|\bU_h^{i}\|_{\bX_1})\|\bv_h\|_{\bX_1}.\label{P4}
\end{align}
The remaining terms are estimated following Lemmas \ref{Intermidiate1} and \ref{Intermidiate2}. Thus, employing the inf-sup condition from Lemma \ref{infsup}, combining all the above estimates in (\ref{P1}) and applying Lemma \ref{Intermidiate1}, we arrive at
\begin{align}
\|r_hp^i-P_h^i\|_{L^2(\Omega_f)} \leq C \|\partial_t \e_f^i\|_{L^2(\Omega_f)} + C(h^{k_1}+h^{k_2}+{\Delta t}^{1/2}). 
\end{align} 
Since $\|p^i-P_h^i\|_{L^2(\Omega_f)}\leq\|p^i-r_hp^i\|_{L^2(\Omega_f)}+ \|r_hp^i-P_h^i\|_{L^2(\Omega_f)} $, using approximation property (\ref{projection1}) and Lemma \ref{Intermidiate2}, we achieve the desired result.
\end{proof}
\begin{Remark}
If we further assume $\bu_{tt}\in {L^{\infty}(0,T,L^2(\Omega_f))}$, then we get an improved pressure error estimate:
\begin{align}
\|p^i-P_h^i\|_{L^2(\Omega_f)} \leq C(h^{k_1}+h^{k_2}+{\Delta t} ).
\end{align}
\end{Remark}
\noindent To find the above estimate, we rewrite (\ref{P1}) as
\begin{align}
(\partial_t \bLambda^i,\bv_h)+& b_{NS}(\bv_h,r_hp^i-P_h^i) =  -(\bu_t^i -\partial_t \bu^i ,\bv_h) - (\partial_t (\bu^i-\bP_h\bu^i),\bv_h) - a_{NS}( \bu^i-\bP_h \bu^i,\bv_h)\nonumber\\
 &- a_{NS}(\bLambda^i,\bv_h)  - \frac{1}{G}(\e_f^i\cdot\btau_{\gamma},\bv_h\cdot\btau_{\gamma})_{{\Gamma}} -  g(e_p^i,\bv_h\cdot\bn_{\gamma})_{{\Gamma}} - c_{NS}(\bu^i,  \bu^i;\bu^i, \bv_h)\nonumber\\
 &
+ c_{NS}(\bU_h^{i-1},  \bU_h^{i-1};\bU_h^{i},\bv_h ) -  b_{NS}(\bv_h,p^i-r_hp^i).\label{P5}
\end{align}
\noindent The term $(\bu_t^i -\partial_t \bu^i ,\bv_h)$ is bounded by applying (\ref{T8}), (\ref{Taylor}) and the Cauchy-Schwarz inequality as follows
\begin{align}
\Big|(\bu_t^i -\partial_t \bu^i ,\bv_h)\Big|\leq C {\Delta}t^{1/2} \bigg(\int_{t_{n-1}}^{t_{n}} |\bu_{tt}|^2_{L^2(\Omega_f)} dt \bigg)^{1/2}\|\bv_h\|_{\bX_1}
\leq C\Delta t \sup_{1\leq i\leq N} \|\bu^i_{tt}\|_{L^2(\Omega_f)}\|\bv_h\|_{\bX_1}.
\end{align}
\noindent For the term $(\partial_t (\bu^i-\bP_h\bu^i),\bv_h)$, using (\ref{T8}), (\ref{Poincare}), and approximation property (\ref{projection2}), we obtain
\begin{align}
\Big|(\partial_t (\bu^i-\bP_h\bu^i),\bv_h)\Big|\leq \frac{1}{\Delta t}\int_{t_{n-1}}^{t_{n}}\|\bu_t-\bP_h\bu_t\|_{L^2(\Omega_f)}\|\bv_h\|_{\bX_1}dt 
\leq Ch^{k_1}\sup_{1\leq i\leq N}|\bu_t^i|_{{k_1+1},\Omega_f}\|\bv_h\|_{\bX_1}.
\end{align}
\noindent The remaining terms on the right-hand side of (\ref{P5}) can be estimated similarly to Theorem \ref{Intermidiate2}. Now Lemma \ref{infsup} yields
\begin{align}
\|\partial_t \bLambda^i\|_{L^2(\Omega_f)}+\gamma^{*}\|r_hp^i-P_h^i\|_{L^2(\Omega_f)} \leq \sup_{\bv_h\in\bX_h\setminus \{0\}} \frac{(\partial_t \bLambda^i,\bv_h)}{\|\bv_h\|_{\bX_1}} + \sup_{\bv_h\in\bX_1\setminus \{0\}} \frac{(\bv_h,r_hp^i-P_h^i)}{\|\bv_h\|_{\bX_1}} .
\end{align}
\noindent The desired estimate is obtained by applying the above bounds in (\ref{P5}).
\begin{Remark}\label{3D}
We can easily extend the optimal results discussed in this study for the three dimensional SIPG case under some certain conditions. The analysis would be modified according to trace inequalities \cite{KJ98, DDE11}, Sobolev embeddings \cite{LS03}, inverse inequalities \cite{DDE11}, etc. Mainly, the estimates of trilinear terms will be changed. To estimate the $L^6$-norm, the existing $L^\mathfrak{p}$ bounds for the three-dimensional case \cite{LS03} can be used:
\begin{align}
\|\bv\|_{L^\mathfrak{p}(\Omega_f)}\leq C_\mathfrak{p}\|\bv\|_{\bX_1}, \,\,\, \forall\bv\in\bX_1, \,\,\, \mathfrak{p}\in [2,6].\label{3D1}
\end{align}
\noindent For the $L^3$-norm estimate, we need a modified bound for the elements in the discrete discontinuous space $\bX_h$ (\cite[eq. (25)]{KHR23}):
\begin{align*}
\|\bv_h\|_{L^3(\Omega)}\leq C \|\bv_h\|^{1/2}_{L^2(\Omega)}\|\bv_h\|^{1/2}_{\bX_1},\quad \forall\bv_h\in\bX_h .
\end{align*}
\noindent In the two-dimensional case, we use an improved $L^4$-norm bound (Lemma $3.8$). But in the three-dimensional case, we have to rely on the $L^3$ and $L^6$-norm bounds. So, using the aforementioned estimate, we can modify the lemma for the $L^3$-norm. For the $L^6$-norm, there is only an estimate (\ref{3D1}). Consequently, we can estimate the nonlinear terms accordingly and derive our desired estimates in the three-dimensional case.  
\end{Remark}
\section{Numerical Results}\label{s7.}
In this section, numerical experiments are provided to validate our theoretical findings. For space discretization, $\mathbb{P}_1-\mathbb{P}_0$ and $\mathbb{P}_2-\mathbb{P}_1$ discontinuous finite element spaces have been employed in the fluid region, and for porous media, $\mathbb{P}_1$ and $\mathbb{P}_2$ discontinuous finite element spaces have been applied, respectively. Time discretization has been done employing the backward Euler scheme. In the first two examples, we take all the values of the parameters as $1$. These examples are computed on the time interval $[0, 1]$ with the final time $T = 1$. Also, we have provided a model test in the last example.\par
\noindent\textbf{Example 6.1} Suppose that the computational domain $\Omega$ is composed of two domains $\Omega_f=[0,1]\times[1,2]$ and $\Omega_p = [0,1]\times [0,1]$. The interface is $\Gamma=(0,1)\times \{1\}$. The exact solution is given by
\begin{align*}
     \bu &= [(x^2(y-1)^2+y)cos(t),-\frac{2}{3}x(y-1)^3+(2-\pi sin(\pi x))cos(t)],  \\
     p &=(2-\pi sin(\pi x))sin(0.5\pi y)cos(t), \quad
     \phi = (2-\pi sin(\pi x))(1-y-cos(\pi y))cos(t).
\end{align*}
\noindent The initial conditions and the source terms are derived from the exact solution. \par
\begin{table}[ht!]
		\centering
		\caption{Numerical results for velocity and pressure with $\mathbb{P}_1$--$\mathbb{P}_0$ element of Example 6.1.}
		\vspace*{-2mm}
		\begin{tabular}{|c|c c| c c| c c| c c|}
		
			\hline

$h$ & $\|\bu(T)-\bU^n\|_{\bX_1}$ & Rate & $\|\bu(T)-\bU^n\|_{L^2(\Omega_f)}$ & Rate & $\|p(T)-P^n\|_{L^2(\Omega_f)}$ & Rate \\[1ex]
			
	\hline   1/2  & $7.1648e-01$  &    $-$      & $1.0965e-01$ &  $-$   & $1.175e-0$ &      $-$       \\
			 1/4  & $4.6875e-01$  &  0.612  & $3.7037e-02$ & 1.566 & $2.5362e-01$ &  1.012  \\
			 1/8  & $1.8518e-01$  &  1.339  & $7.5758e-03$ & 2.289 & $3.1149e-01$ & $-$  \\
			 1/16 & $1.0158e-01$  &  0.866  & $2.0612e-03$ & 1.877 & $2.1920e-01$ & 0.507   \\
			 1/32 & $4.9478e-02$  &  1.018  & $5.0894e-04$ & 2.017 & $1.1302e-01$ & 0.955    \\
			 1/64 & $2.4429e-02$  &  1.019  & $1.2527e-04$ & 2.022 & $5.5286e-02$ & 1.032  \\
			\hline

		\end{tabular}

	\end{table}
	\begin{table}[ht!]
		\centering
		\caption{Numerical results for piezometric head with $\mathbb{P}_1$ element for Example 6.1.}
		\vspace*{-2mm}
		\begin{tabular}{|c|c c| c c| c c|}
		
			\hline

$h$ & $\|\phi(T)-\Phi^n\|_{M_{1p}}$ & Rate & $\|\phi(T)-\Phi^n\|_{L^2(\Omega_p)}$ & Rate \\[1ex]
			
	\hline   1/2  & $4.7204e-01$  &    $-$      & $1.0852e-01$ &             $-$ \\
			 1/4  & $2.8445e-01$ & 0.730  & $3.0102e-02$ & 1.850 \\
			 1/8  & $1.1012e-01$ & 1.368  & $6.7024e-03$ & 2.171 \\
			 1/16 & $5.1523e-02$ & 1.020  & $1.6093e-03$ & 2.032  \\
			 1/32 & $2.4802e-02$ & 1.054  & $4.1222e-04$ & 2.054  \\
			 1/64 & $1.2236e-02$ & 1.095  & $1.0069e-04$ & 2.073   \\
			\hline

		\end{tabular}

	\end{table}
	\begin{table}[ht!]
		\centering
		\caption{Numerical results for velocity  and pressure  with $\mathbb{P}_2$--$\mathbb{P}_1$ element for Example 6.1.}
		\vspace*{-2mm}
		\begin{tabular}{|c|c c| c c| c c|}
		
			\hline

			 $h$ & $\|\bu(T)-\bU^n\|_{\bX_1}$ & Rate & $\|\bu(T)-\bU^n\|_{L^2(\Omega_f)}$ & Rate &  $\|p(T)-P^n\|_{L^2(\Omega_f)}$ & Rate\\[1ex]
			
	\hline   1/2  & $1.3880e-01$  &   $-$       & $1.7166e{-02}$ &          $-$ & $2.9717e-01$  & $-$ \\
			 1/4  & $4.2303e-02$  &  1.714  & $2.7940e{-03}$ & 2.619  & $8.4598e-02$ & 1.813 \\
			 1/8  & $1.3400e-02$  &  1.658  & $3.9836e{-04}$ & 2.810 & $2.826e-02$ & 1.582 \\
			 1/16 & $2.6350e-03$  &  2.046  & $5.0531e{-05}$ & 2.778  & $7.0252e-03$ & 2.008 \\
			 1/32 & $6.3640e-04$  &  2.049  & $7.132e{-06}$ & 2.825 & $1.6817e-03$ & 2.062 \\
			\hline

		\end{tabular}

	\end{table}
	\begin{table}[ht!]
		\centering
		\caption{Numerical results for piezometric head with $\mathbb{P}_2$ element for Example 6.1.}
		\vspace*{-2mm}
		\begin{tabular}{|c|c c| c c| c c|}
		
			\hline

$h$ & $\|\phi(T)-\Phi^n\|_{M_{1p}}$ & Rate & $\|\phi(T)-\Phi^n\|_{L^2(\Omega_p)}$ & Rate \\[1ex]
			
	\hline   1/2  & $2.7202e-01$  &   $-$       & $2.3120e-02$ &             $-$ \\
			 1/4  & $5.7163e-02$ & 2.025  & $2.2967e-03$ & 3.33 \\
			 1/8  & $8.8259e-03$ & 2.069  & $1.8152e-04$ & 3.16 \\
			 1/16 & $1.7428e-03$ & 2.041  & $1.8999e-05$ & 3.05  \\
			 1/32 & $4.222e-04$ & 2.102  & $2.3903e-06$ & 3.25  \\
			\hline

		\end{tabular}

		\vspace{.1cm}

	\end{table}	
	\noindent Tables $1$-$4$ present the computational errors and order of convergence for discontinuous $\mathbb{P}_1$--$\mathbb{P}_0$--$\mathbb{P}_1$ and $\mathbb{P}_2$--$\mathbb{P}_1$--$\mathbb{P}_2$ spaces with the penalty parameters $\sigma_e=10$ and $\sigma_e=40$, respectively. The time step $\Delta t$ has been taken as $O(h)$ and $O(h^2)$ for $\mathbb{P}_1$--$\mathbb{P}_0$--$\mathbb{P}_1$ and $\mathbb{P}_2$--$\mathbb{P}_1$--$\mathbb{P}_2$, respectively. We observe that the numerical results of Tables $1$-$4$ support our theoretically obtained convergence rates, which are of O($h^{k+1}$) for velocity and piezometric head and O($h^{k}$) for pressure in the $L^\infty(L^2)$-norm, respectively.\par 
\noindent\textbf{Example 6.2} The computational domain and the interface condition are taken similarly to Example 6.1. The exact solution is given by
\begin{align*}
 \bu &= ((exp(x+y)+y)cos(t),-(exp(x+y)+x)cos(t)),\\ p &=((cos(\pi x)cos(\pi y)+x)cos(t),\,\,\, \phi = (exp(x+y) - cos(\pi x) + xy)cos(t).
 \end{align*}
 \vspace*{-8mm}
 \begin{table}[ht!]
		\centering
		\caption{Numerical results for velocity and pressure with $\mathbb{P}_1$--$\mathbb{P}_0$ element for Example 6.2.}
		\vspace{-2mm}

		\begin{tabular}{|c|c c| c c| c c| }
		
			\hline
			
			 $h$ & $\|\bu(T)-\bU^n\|_{\bX_1}$ & Rate & $\|\bu(T)-\bU^n\|_{L^2(\Omega_f)}$ & Rate & $\|p(T)-P^n\|_{L^2(\Omega_f)}$ & Rate\\[1ex]
			
	\hline   1/2  & $1.1441e-0$  &   $-$       & $1.6255e-01$ &  $-$    & $2.0606e-0$  &       $-$           \\
			 1/4  & $4.221e-01$  &  1.043  & $3.7043e-02$ & 2.133 & $7.6770e-01$ & 1.424 \\
			 1/8  & $1.9131e-01$  &  1.141  & $7.3448e-03$ & 2.334 & $3.1689e-01$ & 1.276  \\
			 1/16 & $8.2789e-02$  &  1.020  & $1.7469e-03$ & 2.071 & $1.4123e-01$ & 1.065  \\
			 1/32 & $4.0239e-02$  & 1.034  & $4.5824e-04$ & 2.080 & $7.5768e-02$ & 1.083  \\
			\hline

		\end{tabular}

	\end{table}
	\begin{table}[ht!]
		\centering
		\caption{Numerical results for piezometric head with $\mathbb{P}_1$ element for Example 6.2.}
		\vspace{-2mm}
		
		\begin{tabular}{|c|c c| c c| c c|}
		
			\hline
			
			 $h$ & $\|\phi(T)-\Phi^n\|_{M_{1p}}$ & Rate & $\|\phi(T)-\Phi^n\|_{L^2(\Omega_p)}$ & Rate \\[1ex]
			
	\hline   1/2  & $2.9704e-01$  &   $-$       & $6.0629e-02$ &             $-$ \\
			 1/4  & $1.6189e-01$ & 0.875  & $1.1497e-02$ & 2.398 \\
			 1/8  & $5.8120e-02$ & 1.318  & $2.2445e-03$ & 2.356 \\
			 1/16 & $2.3305e-02$ & 1.015  & $4.2782e-04$ & 2.091  \\
			 1/32 & $1.1449e-02$ & 1.028  & $1.0493e-04$ & 2.104  \\
			\hline

		\end{tabular}

		\vspace{.1cm}

	\end{table}
	\begin{table}[ht!]
		\centering
		\caption{Numerical results for velocity and pressure with $\mathbb{P}_2$--$\mathbb{P}_1$ element for Example 6.2.}
		\vspace{-2mm}
		
		\begin{tabular}{|c|c c| c c| c c|}
		
			\hline
			
			 $h$ & $\|\bu(T)-\bU^n\|_{\bX_1}$ & Rate & $\|\bu(T)-\bU^n\|_{L^2(\Omega_f)}$ & Rate & $\|p(T)-P^n\|_{L^2(\Omega_f)}$ & Rate\\[1ex]
			
	\hline   1/2  & $1.8430e-01$  &    $-$      & $2.0923e-02$ &  $-$    & $1.9630e-0$  &    $-$   \\
			 1/4  & $2.9855e-02$  &  2.626  & $2.8857e-03$ & 2.858 & $2.6087e-01$ & 1.912 \\
			 1/8  & $5.4106e-03$  &  2.464  & $3.6913e-04$ & 2.966 & $3.4555e-02$ & 1.741  \\
			 1/16 & $1.0400e-03$  &  2.022  & $4.5288e-05$ & 2.783  & $4.8584e-03$ & 1.833  \\
			 1/32 & $2.4660e-04$  &  2.076  & $5.7404e-06$ & 2.979 & $8.3288e-04$ & 1.982  \\
			\hline

		\end{tabular}

	\end{table}
	\begin{table}[ht!]
		\centering
		\caption{Numerical results for piezometric head with $\mathbb{P}_2$ element for Example 6.2.}
		\vspace*{-2mm}
		\begin{tabular}{|c|c c| c c| c c|}
		
			\hline
			
			 $h$ & $\|\phi(T)-\Phi^n\|_{M_{1p}}$ & Rate & $\|\phi(T)-\Phi^n\|_{L^2(\Omega_p)}$ & Rate \\[1ex]
			
	\hline   1/2  & $2.5102e-01$  &   $-$       & $1.8760e-02$ &             $-$ \\
			 1/4  & $4.8143e-02$ & 2.054  & $2.2560e-03$ & 3.05 \\
			 1/8  & $7.6518e-03$ & 2.091  & $2.7870e-04$ & 3.01 \\
			 1/16 & $2.4781e-03$ & 2.103  & $3.4990e-05$ & 2.98 \\
			 1/32 & $3.1046e-04$ & 2.108  & $4.4059e-06$ & 2.99  \\
			\hline

		\end{tabular}

		\vspace{.1cm}

	\end{table}
	
\noindent Tables $5$-$8$ highlight the computational errors and order of convergence for discontinuous $\mathbb{P}_1$--$\mathbb{P}_0$--$\mathbb{P}_1$ and $\mathbb{P}_2$--$\mathbb{P}_1$--$\mathbb{P}_2$ spaces with the penalty parameters $\sigma_e=10$ and $\sigma_e=40$, respectively. The time step has been chosen similarly to Example 6.1.  
We observe that these results are consistent with the theoretical expectations outlined in Theorems \ref{velo1}, \ref{pressurefinal}, and Lemma \ref{Intermidiate1}.
 \begin{figure}[ht!]
  \centering
  \begin{subfigure}[b]{0.48\textwidth}
    \centering
    \includegraphics[width=\textwidth]{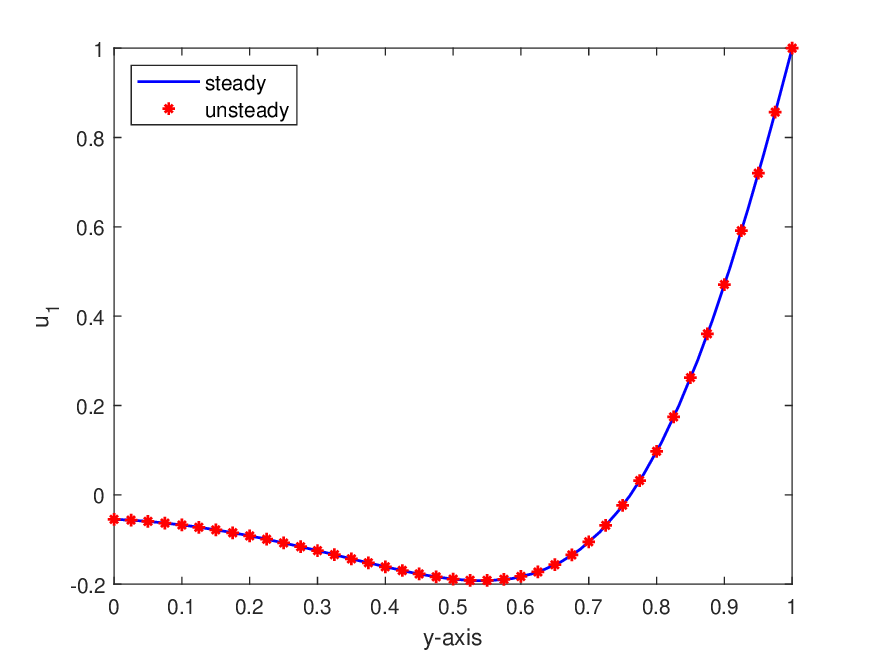}
    \caption{First velocity component along the line $x=0.5$.}
  \end{subfigure}
 \hfill
  \begin{subfigure}[b]{0.48\textwidth}
    \centering
    \includegraphics[width=\textwidth]{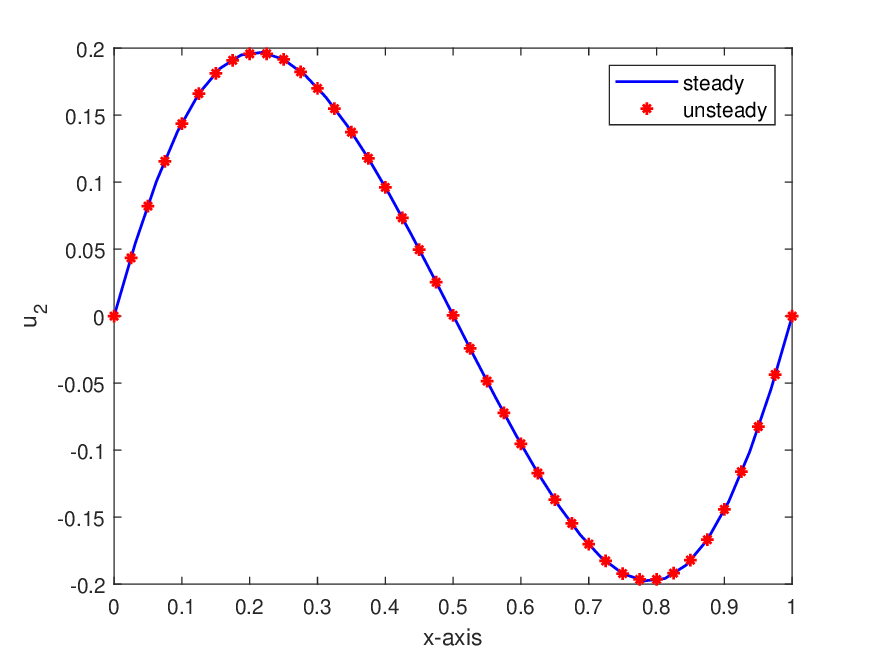}
    \caption{Second velocity component along the line $y=1.5$.}
    \label{fig:second}
  \end{subfigure}
  \centering
  \begin{subfigure}[b]{0.48\textwidth}
    \centering
    \includegraphics[width=\textwidth]{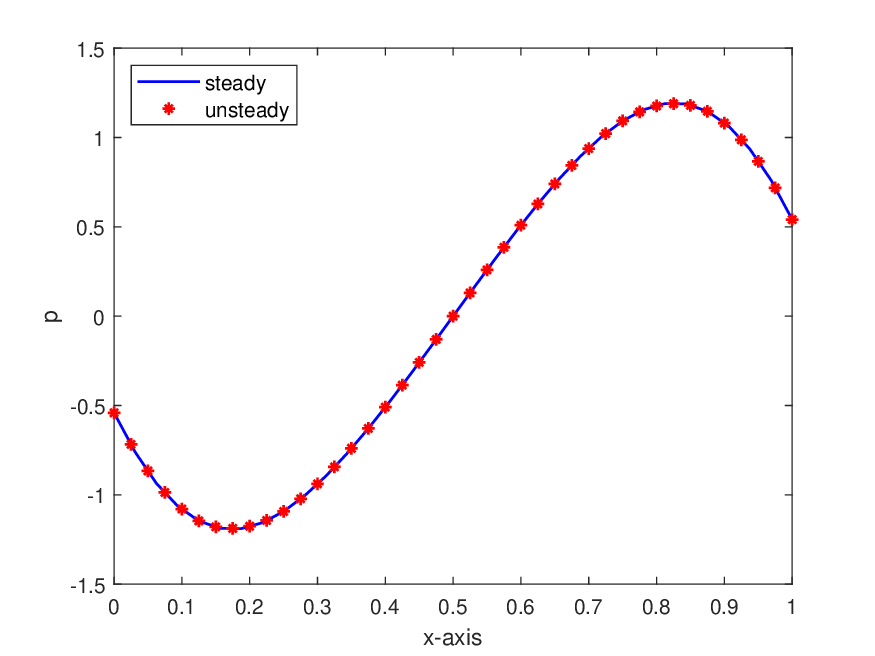}
    \caption{Pressure along the line $y=1.5$.}
  \end{subfigure}
 \hfill
  \begin{subfigure}[b]{0.48\textwidth}
    \centering
    \includegraphics[width=\textwidth]{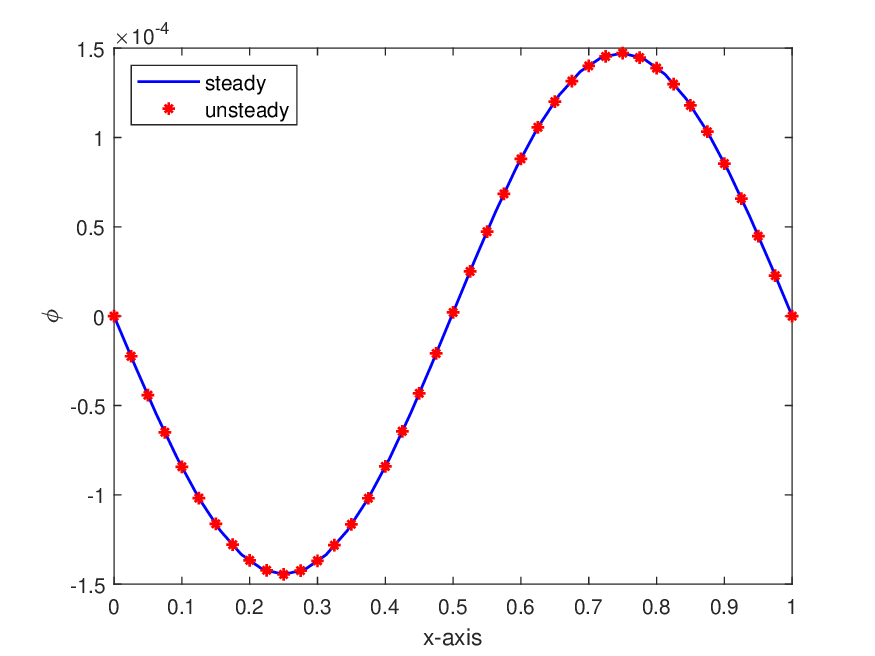}
    \caption{Piezometric head along the line $y=0.5$.}
    \label{fig:second1}
  \end{subfigure}
  \caption{Velocity components, pressure and piezometric head for lid-driven cavity flow with $\mu=1$.}
  \label{fig:all}
\end{figure}

\noindent\textbf{Example 6.3 (Benchmark Problem)} In this example, the lid-driven cavity flow test is considered on the computational regions $\Omega_f=[0,1]\times[1,2],$ $\Omega_p = [0,1]\times  [0,1] $ and $\Gamma=[0,1]\times \{1\} $ with $\f_f=(0,0)$ and $f_p=0$. The velocity on the upper boundary of the domain is $\bu=(1,0)$, and no-slip boundary conditions are considered on the other part of the domain. 
Below, the mentioned boundary conditions are presented:
\begin{align*}
      \begin{cases}
     \bu = (1,0) & \text{on}\,\,[0,1]\times  \{2\},\\
     \bu =(0,0) & \text{on}\,\, \{0\}\times[1,2]\cup \{1\}\times[1,2], \\
     \phi=0 & \text{on}\,\, \partial\Omega_p\setminus\Gamma.
      \end{cases}
\end{align*}
For the simulation, discontinuous $\mathbb{P}_2-\mathbb{P}_1-\mathbb{P}_2$ finite element space has been taken with final time $T=50$, $\Delta t=0.01$, and $h=\frac{1}{32}$. Here we have considered different $\mu=\{1,1/100,1/500,1/1000\}$, $\sigma_e=40$, and other parameters are $1$. In Figure $2$, we have shown a comparison along the lines $(0.5, y)$, $(x, 1.5)$ and $(x, 0.5)$ between unsteady and steady NSD problems. We observe that for large time, the unsteady NSD approximate solution converges to the steady NSD approximate solution. Figure $3$ gives the streamlines for the time-dependent problem for various values of $\mu$. As values of $\mu$ decrease, we notice that the primary vortex moves toward the cavity’s center. Also, a second and third vortex may also develop near the interface.
\begin{figure}[!ht]
  \centering
  \begin{subfigure}[b]{0.48\textwidth}
    \centering
\includegraphics[width=\textwidth]{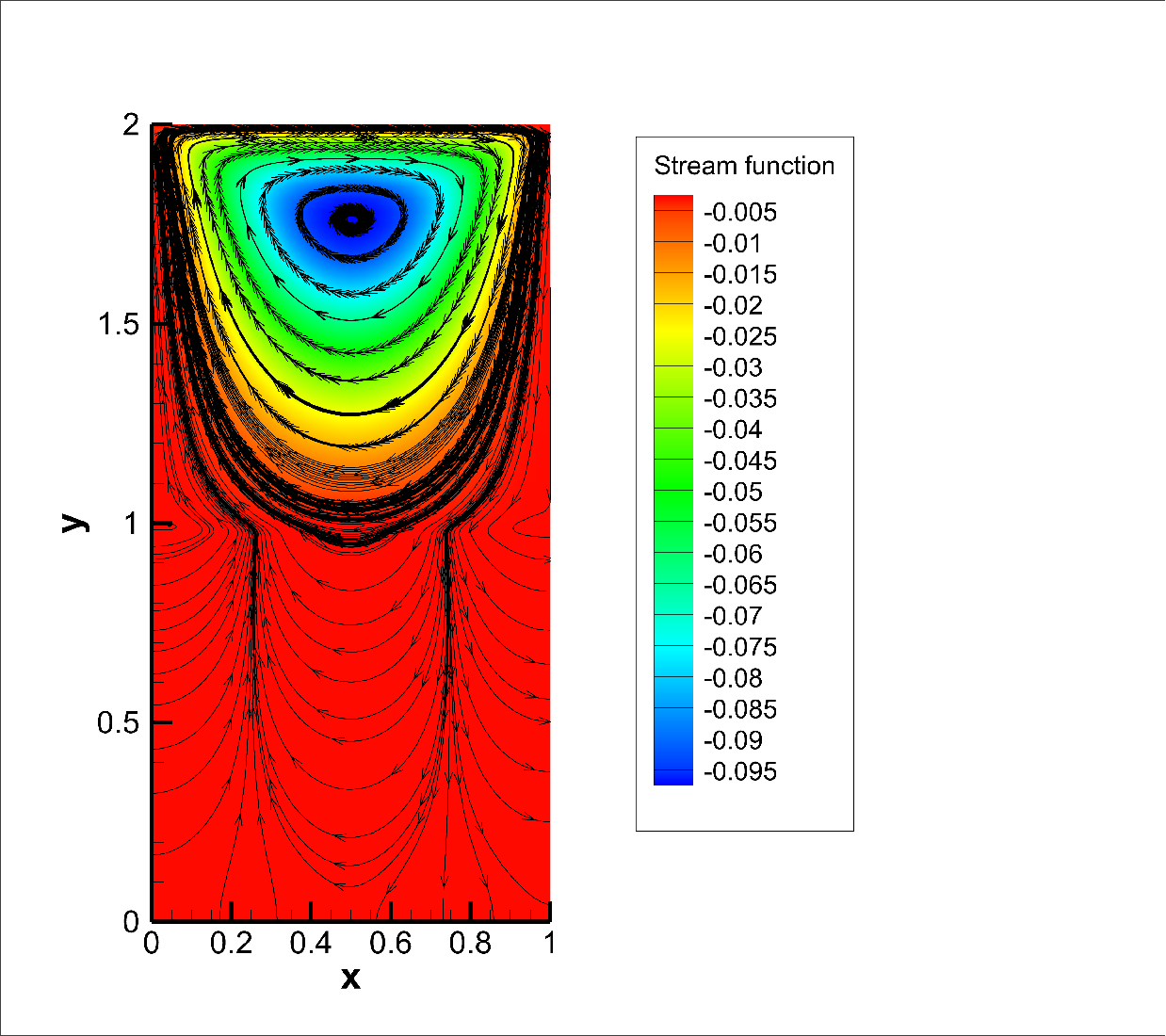}
\caption{$\mu = 1$}\label{thelabel}
\end{subfigure}
\hfill
  \begin{subfigure}[b]{0.48\textwidth}
    \centering
\includegraphics[width=\textwidth]{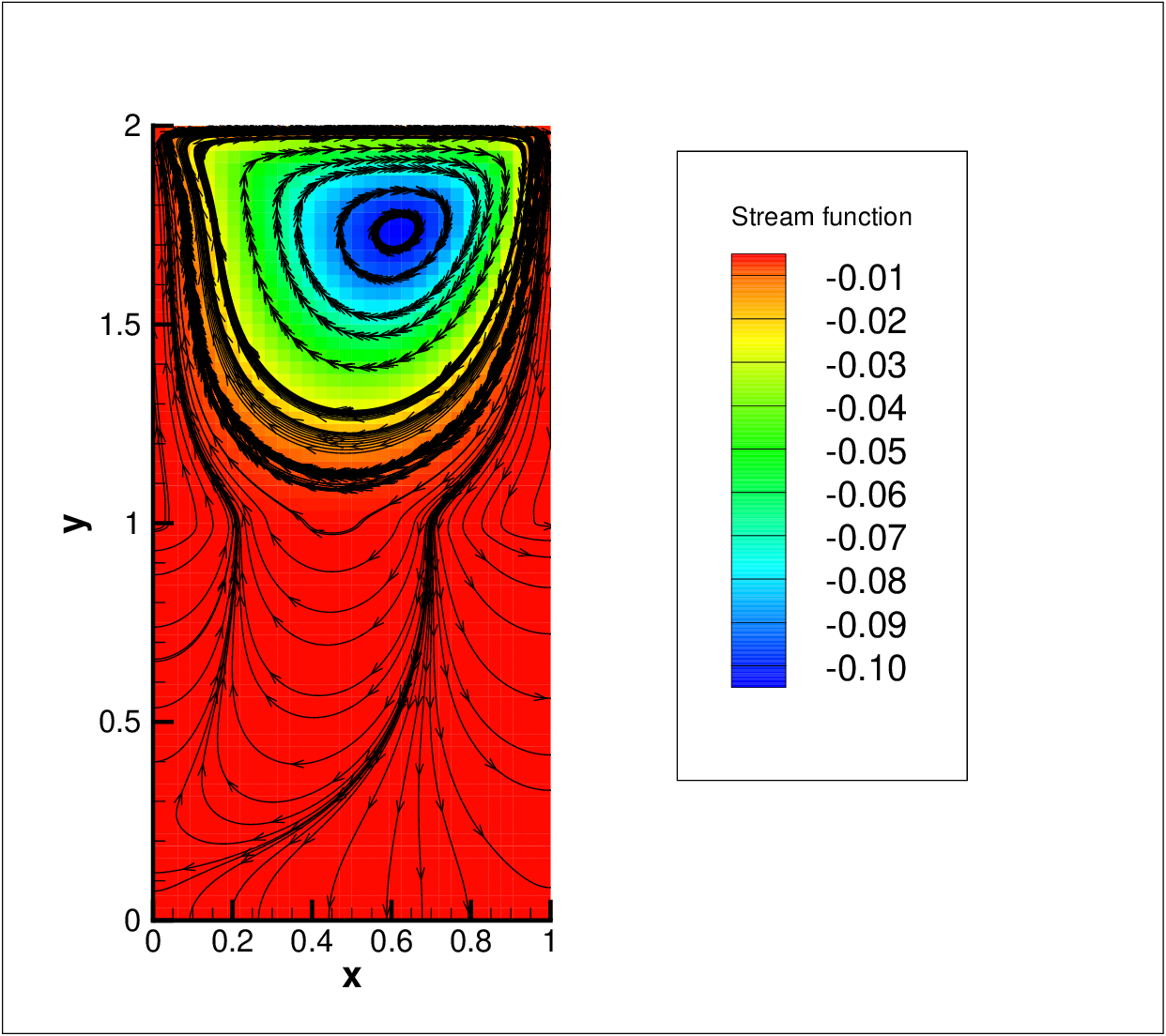} 
\caption{$\mu = 1/100$}\label{thelabel1}
\end{subfigure}
\hfill
  \begin{subfigure}[b]{0.48\textwidth}
    \centering
\includegraphics[width=\textwidth]{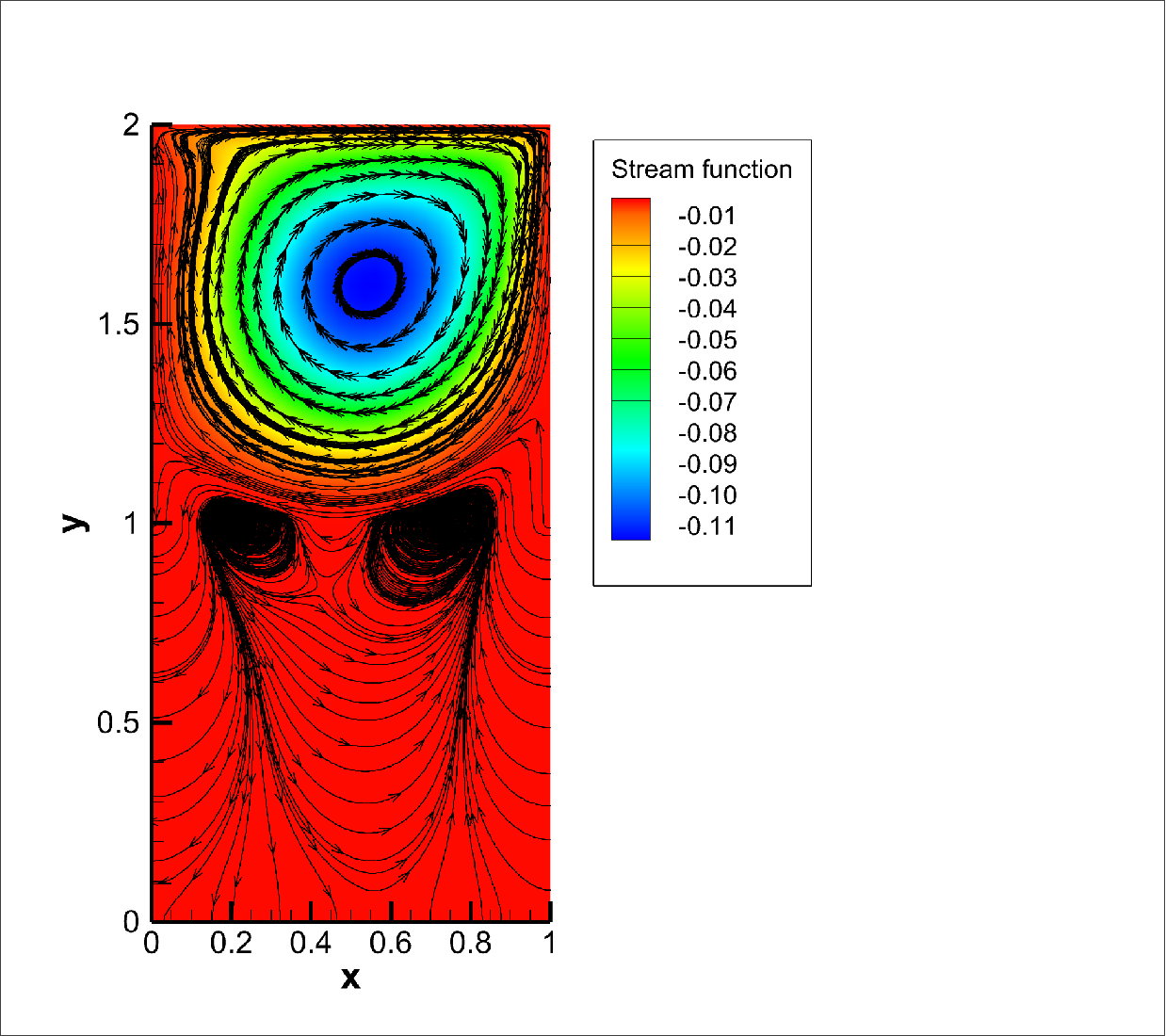} 
\caption{$\mu = 1/500$}\label{thelabel3}
\end{subfigure}
\hfill
  \begin{subfigure}[b]{0.48\textwidth}
    \centering
\includegraphics[width=\textwidth]{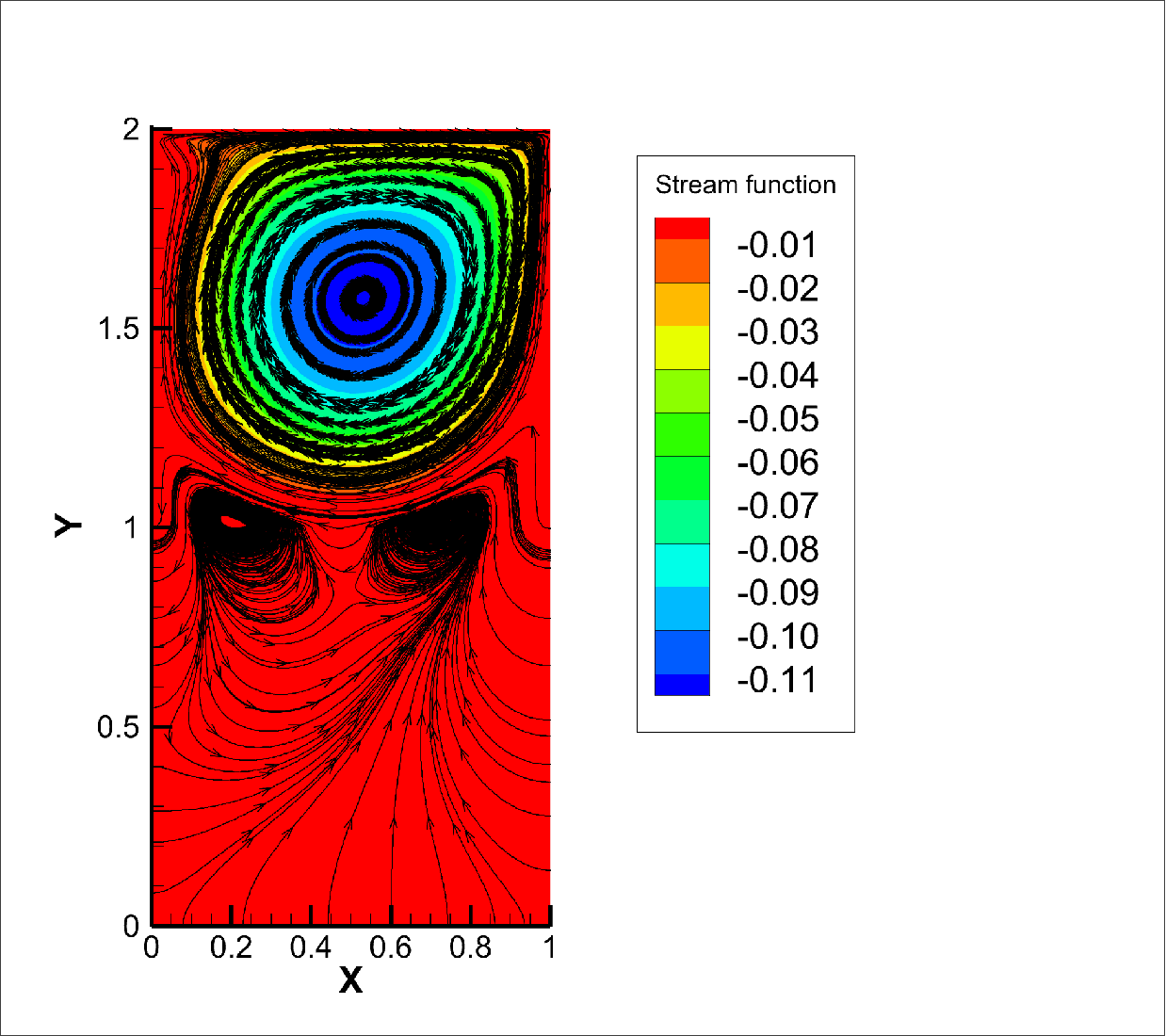} 
\caption{$\mu = 1/1000$}\label{thelabel4}
\end{subfigure}
\caption{Stream lines for evolutionary NSD with $\mu=1,1/100,1/500,1/1000$ at final time $T=50$.}
  \label{fig:all1}
\end{figure}

\vspace*{-2mm}
\section{Conclusion}\label{s8.}
This article analyses a fully discrete DG finite element method for the evolutionary NSD coupled problem. Since standard projections do not yield the optimal results in the $L^\infty(L^2)$-norm, we have introduced modified Oseen operators and established their approximation properties. Then, we have derived optimal error bounds for velocity and piezometric head in $L^\infty(\textbf{L}^2)$ and $L^\infty({L}^2)$-norms as well as energy norm.
Also, a $L^\infty(L^2)$-norm error bound for pressure is obtained. For the temporal discretization, the backward Euler method has been considered.
Finally, we have performed the numerical experiments to confirm our theoretical findings and prove the efficiency of the scheme.
\section*{Acknowledgements} 
The first author would like to express her heartfelt gratitude to the Department of Science and Technology (DST), Government of India, for the financial support (DST/INSPIRE Fellowship/IF210197). 
The second and third authors would like to express their gratitude to Anusandhan National Research Foundation (ANRF), Government of India, under the Core Research Grant project no. CRG/2023/006366.	

\section*{Declarations}
\noindent \textbf{Conflict of interest} The authors declare no competing interests.

\end{document}